\documentclass{article}

\usepackage{times,a4,geometry}
\usepackage{amssymb,amsfonts,amsmath,amscd,amsthm}
\usepackage{multicol,makeidx}
\usepackage{hyperref}

\numberwithin{equation}{subsection}

\newcommand{\dirac}{{\calD}}
\newcommand{\conn}{{\calA(\det P)}}
\newcommand{\SpinC}{\ensuremath{\mathrm{Spin}^{\C}}-structure}
\newcommand{\SpinCs}{{Spin$^{\C}$-structures}}

\newcommand{\R}{\ensuremath{\mathbb{R}}}
\newcommand{\C}{\ensuremath{\mathbb{C}}}
\newcommand{\Z}{\ensuremath{\mathbb{Z}}}
\newcommand{\N}{\ensuremath{\mathbb{N}}}

\newcommand{\CP}{\ensuremath{\mathbb{CP}}}

\newcommand{\Pro}{\ensuremath{\mathbb{P}}}

\newcommand{\sph}{\ensuremath{S}}  
\newcommand{\calC}{\ensuremath{\mathcal{C}}}
\newcommand{\calD}{\ensuremath{\mathcal{D}}}
\newcommand{\calSW}{\ensuremath{\mathcal{SW}}}
\newcommand{\calA}{\ensuremath{\mathcal{A}}}
\newcommand{\calG}{\ensuremath{\mathcal{G}}}
\newcommand{\calB}{\ensuremath{\mathcal{B}}}
\newcommand{\calT}{\ensuremath{\mathcal{T}}}
\newcommand{\calK}{\ensuremath{\mathcal{K}}}
\newcommand{\calH}{\ensuremath{\mathcal{H}}}

\newcommand{\calS}{\ensuremath{\mathcal{S}}}
\newcommand{\calV}{\ensuremath{\mathcal{V}}}
\newcommand{\calE}{\ensuremath{\mathcal{E}}}
\newcommand{\calO}{\ensuremath{\mathcal{O}}}

\newcommand{\M}{\ensuremath{\mathfrak{M}}}
\newcommand{\frakN}{\ensuremath{\mathfrak{N}}}
\newcommand{\preM}{\ensuremath{\widetilde{\mathfrak{M}}}}
\newcommand{\frakJ}{\ensuremath{\mathfrak{J}}}
\newcommand{\frakP}{\ensuremath{\mathfrak{P}}}
\newcommand{\frakX}{\ensuremath{\mathfrak{X}}}
\newcommand{\frakg}{\ensuremath{\mathfrak{g}}}

\newcommand{\im}{\ensuremath{\mathrm{im}}}

\newcommand{\id}{\ensuremath{\mathrm{id}}}
\newcommand{\Id}{\ensuremath{\mathrm{Id}}}
\newcommand{\Ima}{\ensuremath{\mathrm{Im}}}

\newcommand{\vol}{\ensuremath{\mathrm{vol}}}

\newcommand{\End}{\ensuremath{\mathrm{End}}}

\newcommand{\nor}{\ensuremath{\mathrm{nor}}}
\newcommand{\verti}{\ensuremath{\mathrm{vert}}}
\newcommand{\hor}{\ensuremath{\mathrm{hor}}}
\newcommand{\II}{\ensuremath{\mathrm{II}}}

\newcommand{\inj}{\ensuremath{\mathrm{inj}}}
\newcommand{\dist}{\ensuremath{\mathrm{dist}}}
\newcommand{\diam}{\ensuremath{\mathrm{diam}}}
\newcommand{\length}{\ensuremath{\mathrm{L}}}
\newcommand{\Aut}{\ensuremath{\mathrm{Aut}}}
\newcommand{\Lie}{\ensuremath{\mathrm{Lie}}}

\swapnumbers

\theoremstyle{plain} \newtheorem{Def}{\bf D{\footnotesize EFINITION}} [section]
\theoremstyle{plain}\newtheorem{Thm}[Def]{\bf T{\footnotesize HEOREM}} 
\theoremstyle{plain} 
\theoremstyle{plain}\newtheorem{Lem}[Def]{\bf L{\footnotesize EMMA}}
\theoremstyle{plain}\newtheorem{Cor}[Def]{\bf C{\footnotesize
    OROLLARY}} 
\theoremstyle{plain}
\theoremstyle{plain}

\title{\sc On the Riemannian geometry of Seiberg-Witten moduli spaces}
\author{Christian Becker} 
\date{Institut f\"ur Mathematik, Universit\"at Potsdam\\
{\ttfamily becker@math.uni-potsdam.de}}

\begin{document}
\maketitle

\sloppy 

\abstract{
We construct a natural $L^2$-metric on the perturbed Seiberg-Witten moduli
spaces $\M_{\mu^+}$ of a compact 4-manifold $M$, and we
study the resulting Riemannian geometry of $\M_{\mu^+}$. 
We derive a formula which expresses the sectional curvature of $\M_{\mu^+}$ in
terms of the Green operators of the deformation complex of the
Seiberg-Witten equations.
In case $M$ is simply connected, we  construct a Riemannian metric on the Seiberg-Witten
principal $U(1)$ bundle $\frakP \to \M_{\mu^+}$ such that the bundle projection
becomes a Riemannian submersion.
On a K\"ahler surface $M$, the $L^2$-metric on $\M_{\mu^+}$ coincides with the 
natural K\"ahler metric on moduli spaces of vortices.
}

\section{Introduction}
In this paper, we construct a natural $L^2$-metric on the Seiberg-Witten
moduli space $\M_{\mu^+}$ on a compact 4-manifold $M$ with fixed \SpinC{}
$P$. 
The construction follows similar work by {\sc D.~Groisser} and 
{\sc T.~H.~Parker} on the Riemannian geometry of Yang-Mills moduli spaces, see
\cite{groisser-parker87}. 
We study the Riemannian geometry of the $L^2$-metric in the general case of an arbitrary
compact 4-manifold $M$ and in the special case where $M$ is a K\"ahler
surface.  

In the context of Yang-Mills theory, similar research on the geometry of
different (more or less natural) Riemannian metrics on the moduli spaces had
been untertaken in several directions by several people in the works
\cite{groisser-parker87,groisser-parker89,groisser95,babadshanjan-habermann91,habermann92,habermann93,groisser-murray97,maeda-rosenberg-tondeur93,maeda-rosenberg-tondeur95}. 
Although the constructions of the $L^2$-metrics are quite similar, the
naturally arising questions concerning the geometry of the moduli spaces are
rather different for the case Yang-Mills and Seiberg-Witten:
Yang-Mills moduli spaces are noncompact, so it is of particular
interest, whether the natural compactification, which arises
from the analysis of the equations can be realised geometrically, i.e.~as the
completion with respect to the Riemannian distance. One may also ask whether
the volume of the moduli space is finite or infinite, how the metric behaves
near the boundary of the moduli space etc. For results in
these directions (at least for some interesting and accessible examples), we
refer to
\cite{groisser-parker87,groisser-parker89,groisser95,babadshanjan-habermann91,habermann92,habermann93}.
In special cases, where the diffeomorphism type of the moduli space
$\M_{\mu^+}$ can be identified explicitly, one may ask whether the
$L^2$-metric on $\M_{\mu^+}$ coincides with some natural Riemannian metric on
that model space (see e.g.~\cite{groisser-parker87}).  

For the geometry of the Seiberg-Witten moduli space $\M_{\mu^+}$, there arise
different interesting questions: since $\M_{\mu^+}$ is generically a compact
smooth manifold, the $L^2$-metric is always complete and has finite
volume. 
However, since the construction of the moduli space involves the choice of a 
perturbation parameter, one might ask, in how far the $L^2$-metric dependends
on that parameter. 
As the Seiberg-Witten moduli space comes
together with a $U(1)$-bundle $\frakP \to \M_{\mu^+}$, we wonder whether the
construction of the $L^2$-metric on $\M_{\mu^+}$ extends to the total space
$\frakP$. 

On an arbitrary compact smooth 4-manifold, we obtain constructions for quotient 
$L^2$-metrics on the (parametrised) Seiberg-Witten moduli space and (in case $M$ 
is simply connected) on the Seiberg-Witten bundle, which are natural in the sense 
of the following theorem:

\setcounter{section}{3}
\setcounter{Def}{3}

\begin{Thm}
Let $M$ be a compact smooth 4-manifold with a fixed \SpinC{} and
$\mu^+$ (resp.~$\mu^+(t), t \in [0,1]$) generic perturbations such that the
Seiberg-Witten moduli space $\M_{\mu^+}$ (resp.~the parametrised moduli space
$\widehat{\M} = \bigsqcup_{t \in [0,1]} \M_{\mu^+(t)}$) are smooth manifolds
of the expected dimension. 
Then there exists a natural quotient $L^2$-metric on $\M_{\mu^+}$ and a
compatible quotient $L^2$-metric on $\widehat{\M}$ such that the metric
induced from the inclusion of a smooth slice $\M_{\mu^+(t_0)} \hookrightarrow
\widehat{\M}$ is the same as the metric constructed on $\M_{\mu^+(t_0)}$ as
the moduli space with perturbation $\mu^+(t_0)$. 
In case $M$ is simply connected, the Seiberg-Witten bundle $\frakP \to
\M_{\mu^+}$ -- i.e.~the isomorphism class of principal $U(1)$ bundles on
$\M_{\mu^+}$ defining the invariants -- admits a
natural geometric representative carrying a quotient $L^2$-metric such that
the projection $\frakP \to \M_{\mu^+}$ is a Riemannian submersion.
The sectional curvature of those metrics is explicitly given in terms of
the Green operators of the deformation complex of the Seiberg-Witten
equations.
\end{Thm}
\setcounter{section}{1}

On a K\"ahler surface $M$, there is a well known identification of Seiberg-Witten monopoles with vortices.
The vortex equations on compact K\"ahler manifolds were first studied by {\sc Bradlow} and {\sc Garc{\'\i}a-Prada}.
They gave detailed  discussions of existence and uniqueness of vortices and identifications of the corresponding moduli spaces.
In Seiberg-Witten theory, those results yield an identification of the moduli space $\M_{\mu^+}$ with a torus fibration over a complex projective space.
It follows from \cite{garcia-prada_1994a,garcia-prada98}, that the Seiberg-Witten moduli space is a K\"ahler quotient of a K\"ahler submanifold of the configuration space.
As a corollary, our $L^2$-metric is a K\"ahler metric.

The article is organised as follows: In section \ref{chap1}, we briefly recall
the construction of the Seiberg-Witten moduli space and the deformation complex 
we use to construct the $L^2$-metrics.
In section \ref{chap2}, we construct $L^2$-metrics on the Seiberg-Witten bundle 
$\frakP \to \M_{\mu^+}$ and on the parametrised moduli space $\widehat{\M}$. 
We compute a formula for the sectional curvature of $\M_{\mu^+}$, and we show 
that the metric on the slices of the parametrised moduli space coincides with 
the metric constructed before. 
In section \ref{chap3}, we briefly recall the well known identification of the moduli space $\M_{\mu^+}$ as a torus fibration over the complex projective space and as a K\"ahler quotient, which follows from the identification of monopoles as vortices.
We discuss the behaviour of the $L^2$-metric under changes of the perturbation $\mu^+$.
Most of the results in that section are due to {\sc Bradlow}, {\sc Garc{\'\i}a-Prada} and {\sc Hitchin}.

\subsubsection*{Acknowledgements}
The article is based on the authors Ph.D.~thesis. 
His thanks go to the supervisors {\sc Christian B\"ar} and {\sc Paul Gauduchon} for their
encouragement and many stimulating and intriguing discussions, further to the
thesis reviewers {\sc Olivier Biquard} and {\sc Lutz Habermann} for their interest. 
The author thanks the referee for drawing his attention to the K\"ahler quotient construction for moduli spaces of Hermitean-Einstein connections, as given e.g.~in \cite{kobayashi_1987}. 
For financial support, the author thanks the SFB 647 {\it Raum Zeit Materie} of the German Research Foundation.


\section{Notations} \label{chap1}
In this section, we briefly review the construction and basic properties of
the Seiberg-Witten moduli moduli spaces, thereby fixing our notation. As there
exist accessible textbook preparations of that material (our main reference
is \cite{nicolaescu00}), we do not refer to original contributions here.

Throughout this paper, let $M$ be a compact, oriented smooth $4$-manifold
together with a fixed \SpinC{} $P \to M$. Note that this involves more than
just the principal $\mathrm{Spin}^{\C}(4)$ bundle $P$, but to simplify
notation, we only denote it by the symbol $P$.
Let $\calA (\det P)$ be the space of all unitary connections of the
determinant line bundle $\det P$, let $\Sigma^+, \Sigma^-$ be the associated
positive resp.~negative spinor bundle and $\End_0(\Sigma^+)$ the bundle of
tracefree endomorphisms of the positive spinor bundle. The positive Dirac
operator associated with a connection $A \in \calA (\det P)$ is denoted by
$\dirac_A: \Gamma(\Sigma^+) \rightarrow \Gamma(\Sigma^-)$.

The (perturbed) Seiberg-Witten equations are the following coupled 
nonlinear elliptic equations on the configuration space 
$\calC := \conn \times \Gamma(\Sigma^+)$:
\begin{eqnarray}
F_A^+ 
&=& \frac{1}{2} q(\psi,\psi) := (\psi \otimes \psi^*)_0 + \mu^+ \label{SW1}\\
\dirac_A (\psi) 
&=& 0  \; . \label{SW2}
\end{eqnarray}
Here, the 2-form $\mu^+ \in \Omega^2_+(M;i\R)$ is a perturbation parameter.
$F_A^+$ denotes the self-dual part of the curvature $F_A$ of the
connection $A$, and $q$ denotes the real bilinear form
\begin{equation*}
\begin{array}{cccc}
q: 
&\Gamma(\Sigma^+) \times \Gamma(\Sigma^+) 
&\rightarrow 
&\Gamma(\End_0(\Sigma^+)) \\
& & & \\
&(\psi,\phi) 
&\mapsto
&(\psi^* \otimes \phi + \phi^* \otimes \psi)_0 \; .
\end{array}
\end{equation*}
The index $(\cdot)_0$ denotes the trace free part, i.e.
\begin{equation*}
q(\psi,\phi) 
= (\psi^* \otimes \phi + \phi^* \otimes \psi)_0
= \psi^* \otimes \phi + \phi^* \otimes \psi 
  - \frac{1}{2} (\langle \psi,\phi \rangle 
  + \langle \phi,\psi \rangle) \cdot \Id_{\Sigma^+} \; .
\end{equation*}
Solutions of the Seiberg-Witten equations are called 
{\it Seiberg-Witten monopoles} or {\it monopoles}
for short. The space of all monopoles for a fixed perturbation $\mu^+$ is
called the {\it Seiberg-Witten premoduli space} and is denoted by
$\preM_{\mu^+}$. As the zero locus of the Seiberg-Witten map
\begin{equation*}
\begin{array}{cccc}
\calSW_{\mu^+}: 
& \conn \times \Gamma(\Sigma^+)
& \rightarrow
& \Omega^2_+(M;i\R) \times  \Gamma(\Sigma^-) \\
& & &  \\
& \left( \begin{array}{c} A \\ \psi \end{array} \right)
& \mapsto
& \left( \begin{array}{c} F_A^+ - \frac{1}{2}q(\psi,\psi) - \mu^+\\ 
  \dirac_A \psi \end{array} \right) \; , 
\end{array}
\end{equation*}
the premoduli space $\preM_{\mu^+}$ is an infinite dimensional Fr\'echet
submanifold of the configuration space $\calC$ (at least for generic
perturbations). 

The {\it gauge group} $\calG = \Aut(\det P) = \Omega^0(M;U(1))$ acts freely on
the irreducible configuration space 
$\calC^* := \conn \times \big( \Gamma(\Sigma^+) - \{0\} \big)$ by
\begin{equation*}
\calG \ni u:(A,\psi) \mapsto \big( A + 2 u^{-1}du, u^{-1}\psi) \; .
\end{equation*}
The Seiberg-Witten equations are gauge invariant, and the 
{\it Seiberg-Witten moduli space} is defined as the quotient 
\begin{equation*}
\M_{\mu^+} := \preM_{\mu^+} / \calG
\end{equation*}
of the space of monopoles by the gauge group action.

Since the premoduli space $\preM$ is the zero locus of the
Seiberg-Witten map $\calSW_{\mu^+}$, its tangent space in a regular point
$(A,\psi)$ is the kernel of the linearisation in $(A,\psi)$ of
$\calSW_{\mu^+}$. The tangent space in $(A,\psi)$ of the gauge orbit through
$(A,\psi)$ is the image of the linearisation in ${\bf 1} \in \calG$ of the
orbit map through $(A,\psi)$. The linearisation in $(A,\psi)$
of the Seiberg-Witten map $\calSW_{\mu^+}$ is given by:
\begin{equation}
\begin{array}{cccc}
\calT_1:
& \Omega^1(M;i\R) \times \Gamma(\Sigma^+)
& \rightarrow
& \Omega^2_+(M;i\R) \times \Gamma(\Sigma^-) \\
& & & \\
& \left( \begin{array}{c} \nu \\ \phi \end{array} \right)
&\mapsto
& \left( \begin{array}{c}
  d^+\nu - q(\psi,\phi) \\
  \frac{1}{2}\nu \cdot\psi + \dirac_A \phi \end{array} \right)
\end{array}
\end{equation}
The linearisation in ${\bf 1} \in \calG$
of the orbit map through $(A,\psi)$ is given by:
 \begin{equation}
\begin{array}{cccc}
\calT_0:
& \Omega^0(M;i\R)
& \rightarrow
& \Omega^1(M;i\R) \times \Gamma(\Sigma^+) \\
& & & \\
& if
& \mapsto
& \left( \begin{array}{c} 2idf \\  -if \cdot \psi \end{array} \right) 
\end{array}
\end{equation}
Both these linearisations depend on a fixed configuration
$(A,\psi)$ - the one where we linearise the map $\calSW_{\mu^+}$
resp.~where the orbit map is based. We will always drop this
dependence in the notation, but one should keep in mind, that all
formulae derived from these linearisations carry this dependence.

The linearisations $\calT_0, \calT_1$ fit together to the 
elliptic complex $\calK_{(A,\psi)}$, called the deformation complex of the
Seiberg-Witten equations:
\begin{equation}\tag*{$\calK_{(A,\psi)}$}  
0 \longrightarrow
  \Omega^0(M;i\R)
  \stackrel{\calT_0}{\longrightarrow}
  \Omega^1(M;i\R) \times \Gamma(\Sigma^+)
  \stackrel{\calT_1}{\longrightarrow}
  \Omega^2_+(M;i\R) \times  \Gamma(\Sigma^-)
  \longrightarrow 
0 
\end{equation}
It is in reference to this complex, that we denote the linearisations
of the orbit map resp.~the Seiberg-Witten map by $\calT_0$
resp.~$\calT_1$.   

The local structure of the moduli space $\M_{\mu^+}$, especially the 
necessary and sufficient conditions for $\M_{\mu^+}$ to be a
smooth manifold, can easily be described in terms of the elliptic complex
$\calK_{(A,\psi)}$: Since the premoduli space $\preM$ is the zero locus
of the Seiberg-Witten map $\calSW_{\mu^+}$, a necessary condition to apply 
an implicit function theorem is the surjectivity of the map $\calT_1$.
On the other hand, the moduli space is nonsingular only if it does not
contain reducible monopole classes, i.e.~if the orbit map resp.~its
linearisation $\calT_0$ is injective. Thus in the above elliptic complex
$\calK_{(A,\psi)}$ there arise two obstructions for the moduli space
$\M_{\mu^+}$ to be a smooth manifold of the expected dimension near
$(A,\psi)$: the kernel of $\calT_0$  
-- or the zeroth cohomology $\calH^0(\calK_{(A,\psi)})$ of the complex
-- as the obstruction for the gauge action to be free or the moduli space to
be nonsingular, and the cokernel of $\calT_1$ 
-- or the second cohomology $\calH^2(\calK_{(A,\psi)})$
of the complex 
-- as the obstruction for the transversality.

The space $\Gamma_g^+ \subset \Omega^2_+(M;i\R)$ of those perturbations which
admit reducible monopoles $(A,0)$ is a codimension $b_2^+$ hyperplane. 
For $\mu^+ \not\in \Gamma^+_g$, the first obstruction space vanishes. 
One can show that both obstruction spaces vanish for generic perturbations $\mu^+ \in \Omega^2_+(M;i\R)$. 
Thus the moduli space $\M_{\mu^+}$ is generically a smooth manifold of dimension 
$d= - \chi(\calK_{(A,\psi)})$.

Since the Seiberg-Witten equations involve self-dual parts of $2$-forms, the
construction of the moduli space $\M_{\mu^+}$ depends not only on
the perturbation parameter $\mu^+$, but also on the choice of a Riemannian
metric $g$ on $M$. Given two pairs $(\mu^+_0,g_0)$ and $(\mu^+_1,g_1)$ of
generic perturbations and Riemannian metrics, one can show, that for a generic
path $t \mapsto (\mu^+_t,g_t)$ joining $(\mu^+_0,g_0)$ and $(\mu^+_1,g_1)$,
the {\it parametrised moduli space}
\begin{equation*}
\widehat{\M} := \bigsqcup_{t \in [0,1]} \M_{\mu^+(t)} \; .
\end{equation*}
is a smooth cobordism between the smooth moduli spaces $\M_{\mu^+_0}$ and 
$\M_{\mu^+_1}$, if the underlying manifold $M$ satisfies $b_2^+(M) > 1$. 
Elements of $\widehat{\M}$ will be
denoted by $\widehat{[A,\psi]}$. Note that not all the fibres $\M_{\mu^+_t}$ of
the parametrised moduli space need to be smooth manifolds. 

In the case of a $4$-manifold $M$ with $b_2^+(M) = 1$, the perturbations
$\mu^+_i, i=0,1$ may lie on different sides of the separating wall
$\Gamma_g^+$ of reducible perturbations. In that case, the path 
$t \mapsto (\mu^+_t,g_t)$ can be decomposed into subpaths leaving either the
metric or the perturbation fixed. Then the path $t \mapsto \mu^+_t$ may be chosen
in such a way, that it crosses the wall $\Gamma_g^+$ only once. All these
choices can be made generic, and we end up with a singular corbordism joining
the smooth moduli spaces $\M_{\mu^+_0}$ and $\M_{\mu^+_1}$. The singularity in
$\widehat{\M}$ is a cone on a complex projective space. This knowledge yields
explicit formulae for the change of the Seiberg-Witten invariant when
crossing the wall, see e.g.~\cite{nicolaescu00}.

The so called {\it Seiberg-Witten bundle} is an isomorphism class of principal
$U(1)$-bundle over the moduli space $\M_{\mu^+}$, represented by the fibration
\begin{equation}
U(1) 
\hookrightarrow \preM_{\mu^+} / \calG_{x_0} 
\twoheadrightarrow \M_{\mu^+} \, .
\end{equation}
Here $\calG_{x_0}$ is the {\it based gauge group}
\begin{equation}
\calG_{x_0} :=  \{ u \in \calG \; | \; u(x_0) = 1 \} 
\end{equation}
and $x_0 \in M$ is an arbitrary base point. It is easy to see, that for
different base points $x_i \in M, i=0,1$, the representations 
$\varrho_{x_i}: \calG \to U(1), u \mapsto u(x_i)$ are homotopic, and thus that
the quotients $\preM_{\mu^+}/\calG_{x_i}$ are isomorphic as
$U(1)$-bundles. However, on an arbitrary compact, connected $4$-manifold
there are no distinguished base points. Thus the construction of this class of
$U(1)$-bundles seems somewhat ungeometric. 
On a simply connected manifold $M$ we will give a more natural construction 
of a $U(1)$-bundle representing the same isomorphism class. 
This new, geometric representation is needed for the construction of an $L^2$-metric on
the Seiberg-Witten bundle. 


\section{The $L^2$-metric on the moduli space} \label{chap2}
Throughout this chapter, we assume $\mu^+$ to be a perturbation which makes the
moduli space $\M_{\mu^+}$ into a smooth manifold. For a perturbation which
gives rise to a singular moduli space, our construction still yields a
Riemannian metric on the regular part $\M_{\mu^+}^*$ of the moduli space.

We construct a natural $L^2$-metric on the Seiberg-Witten moduli
space $\M_{\mu^+}$ induced from the $L^2$-metric on the configuration space
$\calC$. The tangent space $T_{[A,\psi]}\M_{\mu^+}$ can naturally be
identified with the first cohomology of the elliptic complex
$\calK_{(A,\psi)}$, and we use the elliptic splittings of $T_{(A,\psi)}\calC$
to get an $L^2$-metric on the moduli space $\M_{\mu^+}$. We will allways
assume the perturbation $\mu^+$ to be generic, so that the moduli space embeds
smoothly into the space $\calB^* := \calC^*/\calG$ of gauge
equivalence classes of irreducible configurations.

\subsection{The $L^2$-metric on the configuration space} \label{chap2.1}
The configuration space $\calC = \conn \times \Gamma(\Sigma^+)$
is an affine space, thus it carries a natural $L^2$-metric induced
from the $L^2$-metric on its parallel space
$\Omega^1(M;i\R) \times \Gamma(\Sigma^+)$.
This metric is only a weak Riemannian metric on $\calC$
in the sense that the tangent spaces are not complete with respect to
the $L^2$-topology. A priori, it is not clear whether a weak
Riemannian metric admits a Levi-Civita connection, because the Kozsul
formula gives an element in the cotangent space only. However, on an affine
space there is a natural candidate for a connection, defined by the
directional derivatives:

Let $X,Y \in \frakX (\calC)$ be vector fields on the configuration space,
represented by maps 
{$X,Y: \calC \rightarrow \Omega^1(M;i\R) \times \Gamma(\Sigma^+)$}.
Then the covariant derivative of $Y$ in $(A,\psi)$ in the direction
$X_0 := X_{(A,\psi)}$ is defined by:
\begin{equation}
(\nabla_{X_0} Y)_{(A,\psi)}
:= \frac{d}{dt}\Big|_0 Y\big( \, (A,\psi)  + t X_0 \, \big) \; .
\end{equation}
This connection is obviously torsionfree, it preserves the $L^2$-metric and it
is flat, so we may call it the Levi-Civita connection of the affine space
$\calC$ with respect to the natural $L^2$-metric.  

\subsection{The quotient $L^2$-metric on the moduli space} \label{chap2.2}
We construct a Riemannian metric on the Seiberg-Witten moduli
space $\M_{\mu^+}$, inherited from the quotient metric on the space of
gauge equivalence classes of irreducible configurations 
$\calB^*:= \calC^* / \calG$ via the embedding 
$\M_{\mu^+}\hookrightarrow \calB^*$. 
We outline how to compute the sectional curvature of this metric in terms of
the Green operators of the elliptic complex $\calK_{(A,\psi)}$ associated with
a monopole $(A,\psi)$. 
Similar $L^2$-metrics on the Yang-Mills moduli spaces had been studied by 
{\sc Groisser, Habermann, Matsumoto, Matumoto} and {\sc Parker} with several approaches to different special cases. 
The work \cite{groisser-parker87} of {\sc Groisser} and {\sc Parker} gives a
detailed introduction to the construction of $L^2$-metrics on moduli spaces of
monopoles. 

The gauge group $\calG = \Omega^0(M;U(1))$ acts on $\calC$ by 
$u: (A,\psi) \mapsto (A + 2 u^{-1}du,u^{-1}\psi)$, and the induced action on
$T\calC$ is given by $u:(\nu,\phi) \mapsto (\nu,u^{-1}\phi)$. Hence the
$L^2$-metric on $\calC$ is $\calG$-invariant, and the quotient space
$\calB^* := \calC^* / \calG$ of gauge equivalence classes of irreducible
configurations carries a unique (weak) Riemannian metric such that the
projection $\calC^* \rightarrow \calB^*$ is a Riemannian submersion. 
We use an infinite dimensional analogue of the O'Neill formula for
Riemannian submersions to derive a formula for the sectional curvature of this
quotient metric on $\calB^*$. The Gauss equation for the embedding 
$\M_{P,\mu^+} \hookrightarrow \calB^*$ then yields a formula for the sectional
curvature of the Seiberg-Witten moduli space $\M_{\mu^+}$.
$$
\begin{CD}
\preM_{\mu^+} @>>{\quad}> \calC^* \\
@VVV @VV{\scriptstyle{\mathrm{O'Neill}}}V\\
\M_{\mu^+} @>>{\scriptstyle{\mathrm{Gauss}}}> \calC^*/\calG
\end{CD}
$$
Both the O'Neill formula and the Gauss equation involve orthogonal
projections onto subspaces of the tangent space. These are given in
terms of the Green operators of the elliptic complex $\calK_{(A,\psi)}$
associated with a monopole $(A,\psi)$. The tangent space in $(A,\psi)$ of the
premoduli space $\preM$ is the kernel of the linearisation $\calT_1$ in
$(A,\psi)$ of the Seiberg map $\calSW_{\mu^+}$. Correspondingly, the tangent
space in $(A,\psi)$ of the gauge orbit through $(A,\psi)$ is the image of the
linearisation in  $\calT_0$ of the orbit map through $(A,\psi)$. The tangent
space in $[A,\psi]$ of the moduli space $\M_{\mu^+}$ may then be identified
with the intersection of $\ker \calT_1$ and the orthogonal complement of 
$\im \calT_0$. The ellipticity of the complex $\calK_{(A,\psi)}$ yields the
$L^2$-orthogonal splittings  
\begin{eqnarray}
\Omega^0(M;i\R)
&=& \ker \calT_0 \oplus \im \calT_0^*  \label{ellsplitM0} \\
\Omega^1(M;i\R) \times \Gamma(\Sigma^+)
&=& \ker \calT_0^* \oplus \im  \calT_0   \notag  \\
&=& \ker \calT_1   \oplus \im  \calT_1^* \notag  \\
&=& (\ker \calT_0^* \cap \ker \calT_1) \label{ellsplitM1} 
\oplus \im  \calT_0 \oplus \im  \calT_1^* \\
\Omega^2_+(M;i\R) \times  \Gamma(\Sigma^-)
&=& \ker \calT_1^* \oplus \im \calT_1 \label{ellsplitM2} \;.
\end{eqnarray}
All these operators implicitly depend on the configuration $(A,\psi)$
where we do the linearisations, but we will drop this dependence in
the notation.

The adjoint of $\calT_0$ is the operator
\begin{equation*}
\begin{array}{cccc}
\calT_0^*:
& \Omega^1(M;i\R) \times \Omega^0(M;L)
& \rightarrow
& \Omega^0(M;i\R) \\
& \left( \begin{array}{c} \nu \\ \phi \end{array} \right)
& \mapsto
& 2d^*\nu \, + \, i \Ima \langle \psi,\phi \rangle
\end{array}
\end{equation*}
and the adjoint of $\calT_1$ is the operator
\begin{equation}
\begin{array}{cccc}
\calT_1^*:
& \Omega^2_+(M;i\R) \times \Gamma(\Sigma^-)
& \rightarrow
& \Omega^1(M;i\R) \times \Gamma(\Sigma^+) \\
& & & \\
& \left( \begin{array}{c} \mu \\ \xi \end{array} \right)
& \mapsto
& \left( \begin{array}{c}
  d^*\mu + \frac{i}{2} \Ima \langle (\cdot) \cdot \psi,\xi \rangle \\
  \dirac_A \xi - 2 \mu \cdot \psi
  \end{array} \right)
\end{array}
\end{equation}
where $\langle (\cdot) \cdot \psi,\xi \rangle$ denotes the 1-form
$\frakX(M) \ni X \mapsto \langle X \cdot \psi,\xi
\rangle$.

The Laplacians
\begin{equation*}
\begin{array}{rccc}
L_0 = \calT_0^* \circ \calT_0: &
\Omega^0(M;i\R) &
\rightarrow &
\Omega^0(M;i\R) \\
L_1 = \calT_0 \circ \calT_0^* \oplus \calT_1^* \circ \calT_1: &
\Omega^1(M;i\R) \times \Gamma(\Sigma^+) &
\rightarrow &
\Omega^1(M;i\R) \times \Gamma(\Sigma^+) \\
L_2 = \calT_1 \circ \calT_1^*: &
\Omega^2_+(M;i\R) \times \Gamma(\Sigma^-) &
\rightarrow &
\Omega^2_+(M;i\R) \times \Gamma(\Sigma^-)
\end{array}
\end{equation*}
associated with the elliptic complex $\calK_{(A,\psi)}$ are invertible on the complements of their kernels, and
the Green operators $G_j, j=0,1,2$ are defined as the extensions by $0$ of
those inverses:
\begin{equation*}
\begin{array}{rcccl}
G_0: &
\Omega^0(M;i\R) &
\rightarrow &
\Omega^0(M;i\R) \, , &
G_0:= (L_0|_{\im \calT_0^*})^{-1} \oplus 0 \\
G_1: &
\Omega^1(M;i\R) \times \Gamma(\Sigma^+) &
\rightarrow &
\Omega^1(M;i\R) \times \Gamma(\Sigma^+) \, , &
G_1 := (L_1|_{\im \calT_0 \oplus \calT_1^*})^{-1} \oplus 0 \\
G_2: &
\Omega^2_+(M;i\R) \times \Gamma(\Sigma^-) &
\rightarrow  &
\Omega^2_+(M;i\R) \times \Gamma(\Sigma^-) \, , &
G_2:= (L_2|_{\im \calT_1})^{-1} \oplus 0  \\
\end{array}
\end{equation*}
These Green operators are nonlocal elliptic pseudo-differential operators.

Using the splittings \eqref{ellsplitM0} -- \eqref{ellsplitM2},
the orthogonal projectors onto the vertical space 
$\calV_{(A,\psi)}:= \im \calT_0$
resp.~the horizontal space 
$\calH_{(A,\psi)}:= \ker \calT_0^*$
of the gauge action as well as onto the tangent space
$T_{(A,\psi)}\preM = \ker \calT_1$
resp.~the normal space  
$N_{(A,\psi)}\preM = \im \calT_1^*$
of the premoduli space are given as:
\begin{alignat}{3}
\verti_{(A,\psi)}&= \calT_0 \circ G_0 \circ \calT_0^*
&&\quad
&\hor_{(A,\psi)} 
   &= \id_{\Omega^1(M;i\R) \times \Gamma(\Sigma^+)} - \verti_{(A,\psi)} 
\label{vert-hor}  \\
\tan_{(A,\psi)} 
   &= \id_{\Omega^1(M;i\R) \times \Gamma(\Sigma^+)} - \nor_{(A,\psi)}
&&\quad
&\nor_{(A,\psi)} &= \calT_1^* \circ G_2 \circ \calT_1  \label{tan-nor}
\end{alignat}
That these operators are in fact the orthogonal projections is quite
obvious, since e.g.~the operator
$\verti = \calT_0 \circ G_0 \circ \calT_0^*$
is the identity on $\im \calT_0$ and vanishes on the orthogonal complement
$\ker \calT_0^*$, thus it is the orthogonal projection from
$\Omega ^1(M;i\R) \times \Gamma(\Sigma^+)$ to
$\calV_{(A,\psi)} = \im \calT_0$. 

We thus have the following natural $L^2$-orthogonal splitting of 
the linearised (irreducible) configuration space $T_{(A,\psi)}\calC^*$:
\begin{equation} \label{TMsplit}
T_{(A,\psi)}\calC^*
= \Omega^1(M;i\R) \times \Gamma(\Sigma^+)
= \underbrace{(\ker \calT_0^* \cap \ker \calT_1)}_{\cong \, 
   T_{[A,\psi]}\M}
   \oplus \im \calT_0 \oplus \im  \calT_1^*  \; .
\end{equation}
By restriction of the $L^2$-metric from $T_{(A,\psi)}\calC^*$ 
to the orthogonal direct summand 
$T_{[A,\psi]}\M_{\mu^+}\cong \ker \calT_0^* \cap \ker \calT_1$, we get a
natural $L^2$-metric on the moduli space $\M_{\mu^+}$, which we call the 
{\it quotient $L^2$-metric}.
In section \ref{chap2.5} below, we compute a formula for the sectional
curvature of this metric using the O'Neill formula and the
Gauss equation together with the above identifications of the orthogonal
projectors.

\subsection{The quotient $L^2$-metric on the Seiberg-Witten bundle} 
\label{chap2.3}
We construct a natural $L^2$-metric on the total space
$\frakP$ of the Seiberg-Witten bundle in the same way as we did for the
moduli space by only replacing the orthogonal splitting \eqref{TMsplit}: 
here we must split the tangent space of $\frakP$ from $\Omega^1(M;i\R) \times
\Gamma(\Sigma^+)$. 
To identify the tangent space of $\frakP$ in a way that automatically yields an orthogonal splitting, we construct a new representative of the isomorphism class of $U(1)$-bundles 
$\frakP \rightarrow \M_{\mu^+}$. 
This reperesentative is more natural from the point of view
of the geometry of the $L^2$-metric than one induced by the based gauge group
$\calG_{x_0}$. For this construction, we need $M$ to be simply connected. 

In the representation of the Seiberg-Witten bundle $\frakP \rightarrow
\M_{\mu^+}$ as the quotient $\preM / \calG_{x_0} \rightarrow \M_{\mu^+}$ of the premoduli space by the based gauge group, the $U(1)$-action on $\frakP \rightarrow \M_{\mu^+}$ comes from the action of the constant gauge transformations $U(1) \subset \calG$ on $\preM$. 
Thus to construct the natural quotient $L^2$-metric on $\frakP$ we need to split the gauge group
$L^2$-orthogonally into constant and non-constant gauge transformations. 
The splitting provided by the based gauge group $\calG_{x_0}$ is not appropriate, since the Lie algebra of the based gauge group
\begin{equation*}
\Lie \calG_{x_0} 
= \{ if \in \Omega^0(M;i\R) \; | \; f(x_0) = 0 \}
\end{equation*}
is $L^2$-dense in the Lie algebra 
$\Lie \calG = \Omega^0(M;i\R)$
of the full gauge group. Thus it is not topologically complemented with
respect to the $L^2$-topology. 

However, the Lie algebra $\frakg = \Omega^0(M;i\R)$ splits naturally as the
orthogonal direct sum of the constant functions and those functions, which
integrate to 0 with respect to the volume form induced by the fixed Riemannian
metric $g$:  
\begin{equation*}
\frakg
= \Lie \calG 
= \Omega^0(M;i\R) 
= i\R \oplus \left\{ if \in \Omega^0(M;i\R) 
  \; \left| \; 
  \mbox{$\frac{1}{\vol(M)} \int_M if \, dv_g = 0$} \right\} \right. \; .
\end{equation*}
In case $M$ is sipmly connected, this splitting can be realised via a splitting of the gauge group $\calG$ itself:

\begin{Def}
The {\it reduced gauge group} is the subgroup $\calG_{\infty} \subset \calG$
of all gauge transformations $u \in \calG$, which satisfy
\begin{equation*}
\exp \left( \frac{1}{\vol(M)} \int_M \log u \, dv_g  \right) = 1 \; .
\end{equation*}
\end{Def}
By this definition we obtain a topological splitting -- in the sense of Fr\'echet-Lie groups -- of the gauge group as
$\calG =  U(1) \times \calG_{\infty}$. 
The Lie algebra
\begin{equation*}
\frakg_{\infty}
= \Lie \calG_{\infty}
= \left\{ if \in \Omega^0(M;i\R) 
  \mid \mbox{$\frac{1}{\vol(M)} \int_M if \, dv_g = 0$} \right\} \; .
\end{equation*}
is the orthogonal complement of $i\R = \Lie U(1)$ in $\frakg$.
Thus the Lie algebra of the full gauge group $\calG$ splits $L^2$-orthogonally
as:  
\begin{equation*}
\frakg = \Omega^0(M;i\R) = i\R \oplus \frakg_{\infty} \; .
\end{equation*}
The quotient of the premoduli space $\preM$ by the reduced gauge group
$\calG_{\infty}$ yields another natural $U(1)$-bundle over the moduli
space $\M_{\mu^+}$. 
The $\calG_{\infty}$ equivalence class of a monopole $(A,\psi) \in \preM$ will be denoted by
$[A,\psi]_{\infty}$.
\begin{Lem}
The $U(1)$-bundle $\preM / \calG_{\infty} \rightarrow \M_{\mu^+}$
represents the isomorphism class $\frakP \rightarrow \M_{\mu^+}$.
Thus for any $x_0 \in M$, the $U(1)$-bundles 
$\preM / \calG_{\infty} \rightarrow \M_{\mu^+}$ and 
$\preM / \calG_{x_0} \rightarrow \M_{\mu^+}$ are isomorphic.
\end{Lem}
\begin{proof}
We define two representations 
$\varrho_{x_0},\varrho_{\infty}: \calG \rightarrow U(1)$, whose kernels
are the subgroups $\calG_{x_0}$ resp. $\calG_{\infty}$. We show that the
bundles $\preM / \calG_{x_0} \rightarrow \M_{\mu^+}$ resp. 
$\preM / \calG_{\infty} \rightarrow \M_{\mu^+}$ are associated from the 
principal $\calG$ bundle $\preM \rightarrow \M_{\mu^+}$ via the representations
$\varrho_{x_0}, \varrho_{\infty}$. Then a homotopy of representations from
$\varrho_{x_0}$ to $\varrho_{\infty}$ yields a homotopy of the associated
principal bundles. This implies that $\preM / \calG_{x_0}$ and 
$\preM / \calG_{\infty}$ have the same first Chern class and thus are
isomorphic.  

The representations 
$\varrho_{x_0},\varrho_{\infty}: \calG \rightarrow U(1)$ 
are defined by:
\begin{eqnarray*}
\varrho_{x_0}(u)
&:=& 
u(x_0) \\
\varrho_{\infty}(u)
&:=&
\exp \left( \mbox{$ \frac{1}{\vol(M)} \int_M \log u \, dv_g$} \right)
\end{eqnarray*}
Obviously, $\ker \varrho_{x_0} = \calG_{x_0}$, whereas
$\ker \varrho_{\infty} = \calG_{\infty}$.
To show, that the quotients $\preM / \calG_{x_0}$ resp.~$\preM / \calG_\infty$
are principal $U(1)$-bundles associated from the principal $\calG$-bundle 
$\preM \to \M_{\mu^+}$, we consider the map
\begin{equation*}
\begin{array}{ccc}
\preM \times U(1) 
& \rightarrow
& \preM \\
& & \\
\left( \begin{array}{c}
A \\  \psi \\ \lambda \end{array} \right)
& \mapsto
& \left( \begin{array}{c}
  A \\ \lambda^{-1} \cdot \psi \end{array} \right) \; .
\end{array}
\end{equation*}
The full gauge group $\calG$ acts on the $U(1)$-factor of the left hand side
via the representations $\varrho_{x_0}$ resp.~$\varrho_\infty$. 
This map is equivariant with respect to the action of $\calG$ on the left
hand side and of $\calG_{x_0}$ resp.~$\calG_\infty$ on the right hand
side. By taking quotients, it descends to isomorphisms
\begin{equation*}
\begin{array}{cccc}
\preM \times_{\varrho_{x_0}} U(1) 
& \cong 
& \preM / \calG_{x_0}  \\
& & \\
\preM \times_{\varrho_{\infty}} U(1)  
& \cong 
& \preM / \calG_{\infty} \; .
\end{array}
\end{equation*}

To construct a homotopy of $U(1)$-bundles from $\preM \times_{\varrho_{x_0}}
U(1)$ to  $\preM \times_{\varrho_{\infty}} U(1)$ it suffices to construct a
homotopy of representations from $\varrho_{x_0}$ to $\varrho_{\infty}$. 
We define such a homotopy $H$ as follows:
\begin{equation*}
\begin{array}{cccl}
H: 
& \calG \times [0,1] 
& \rightarrow 
& U(1) \\
& & & \\
& (u,t) 
& \mapsto 
& \begin{cases} 
  u(x_0) 
  & \quad t = 0 \\
  \exp \left( \mbox{$\int_M \rho_t \cdot \log u \, dv_g$} \right) 
  & \quad t \in (0,1] \,.
  \end{cases}
\end{array}
\end{equation*}
Here the family $\rho_t$ is a smoothing of the Dirac distribution, such that for
any function $f$ the integral \mbox{$\int_M \rho_t \cdot f \, dv_g$}
converges to $f(x_0)$ when $t$ tends to $0$ and  
$\rho_t \equiv \frac{1}{\vol(M)}$ for $t$ close to $1$.
Taking a nonnegative smooth function $\gamma: \R \rightarrow \R$, constant near $0$ with support in $[-1,1]$ and $\int_{\R^4} \gamma(|x|)dx =1$, and setting 
$2 \epsilon < \min \{1,\inj(M,g)\}$, we define the family $\rho_t: M \rightarrow \R$ for 
$t \in (0,1]$ as:
\begin{equation*}
\rho_t(x) 
:= \begin{cases}
   \frac{1}{t^4} \cdot \gamma\left(\frac{\dist(x_0,x)}{t}\right)
 & \quad t \in (0,\epsilon] \\
   \left( \frac{-t}{\epsilon} + 2 \right) 
     \cdot \frac{1}{t^4} \cdot \gamma\left(\frac{\dist(x_0,x)}{t}\right) \,
   + \left( \frac{t}{\epsilon} -1 \right) \cdot \frac{1}{\vol(M)}
 & \quad t \in [\epsilon,2\epsilon] \\
   \frac{1}{\vol(M)} 
 & \quad t \in [2\epsilon,1]
\end{cases}
\end{equation*}
Here $\dist(x,x_0)$ denotes the Riemannian distance from $x$ to $x_0$.
By construction, the integral $\int_M \rho_t \cdot \log u \, dv_g$ tends to
$u(x_0)$ as $t$ tends to zero, thus the homotopy $H$ as defined above is
continuous in $t$ and satisfies $H_0 = \varrho_{x_0}$. Since 
$\rho_t \equiv \frac{1}{\vol(M)}$ near $t=1$, we also have 
$H_1 = \varrho_{\infty}$. Thus $H$ is a homotopy from $\varrho_{x_0}$ to
$\varrho_{\infty}$ as claimed. By construction, for any $t \in [0,1]$, the map
$H_t:\calG \rightarrow U(1)$ is a representation. 

The homotopy of representations $H: \calG \times [0,1] \rightarrow U(1)$
defines a homotopy   
\begin{equation*} 
\widehat{\frakP} \rightarrow \M_{\mu^+}\times [0,1] \; , \qquad
\widehat{\frakP}_t := \preM \times_{H_t} U(1) 
\end{equation*}
of $U(1)$-bundles over $\M_{\mu^+}$ from 
$\widehat{\frakP}_0 
= \preM \times_{\varrho_{x_0}} U(1) \cong \preM / \calG_{x_0}$  
to
$\widehat{\frakP}_1 
= \preM \times_{\varrho_{\infty}} U(1) \cong \preM / \calG_{\infty}$.
This implies that these bundles have the same first Chern number and are
thus isomorphic. 
\end{proof}

Now we can construct an $L^2$-metric on the total space $\frakP$ in much the
same way as we did for the moduli space $\M_{\mu^+}$: 
We identify the tangent space $T_{[A,\psi]_{\infty}}\frakP$ as the
intersection of the kernel of $\calT_1$ with the orthogonal complement of the
$\calG_{\infty}$-orbit through $(A,\psi)$. 
Then we split the linearized configuration space $\Omega^1(M;i\R) \times \Gamma(\Sigma^+)$ into $T_{[A,\psi]_{\infty}}\frakP$ and its orthogonal complement.

The reduced gauge group $\calG_{\infty}$ is a Fr\'echet-Lie subgroup of the full gauge group $\calG$, and its Lie algebra $\frakg_{\infty}$ is a tame direct summand of $\frakg$. 
Thus we get a slice theorem for the action of $\calG_{\infty}$ on $\calC^*$ in the tame smooth category in the same way as for the Fr\'echet-Lie group $\calG$ (for those slice theorems in the tame smooth category see \cite{abbati_cirelli_mania_1989} and \cite{subramaniam_1989}). 
When $\calS_{(A,\psi)}$ is a local slice for the full gauge group $\calG$, then 
$U(1)\cdot\calS_{(A,\psi)}$ is a local slice for the reduced gauge group $\calG_{\infty}$. 
The linearisation of the orbit map for $\calG_{\infty}$ is the restriction to $\frakg_{\infty}$ of the linearisation $\calT_0$ of the orbit map for $\calG$. 
The linearisations of the orbit map and of the Seiberg-Witten map fit together into the
complex:   
\begin{equation} \tag*{$\calK^{\infty}_{(A,\psi)}$}
   0 \longrightarrow
     \frakg_{\infty}
     \stackrel{\calT_0|_{\frakg_{\infty}}}{\longrightarrow}
     \Omega^1(M;i\R) \times \Gamma(\Sigma^+)
     \stackrel{\calT_1}{\longrightarrow}
     \Omega^2_+(M;i\R) \times  \Gamma(\Sigma^-)
     \longrightarrow 0
\end{equation}
The adjoint $\calT_0|_{\frakg_{\infty}}^*$ of the restriction 
$\calT_0|_{\frakg_{\infty}}$ is the composition of the adjoint
of $\calT_0$ with the orthogonal projection to $\frakg_{\infty}$:
\begin{equation*}
\begin{array}{cc}
\calT_0|_{\frakg_{\infty}}^*:
&\Omega^1(M;i\R) \times \Gamma(\sigma^+) \rightarrow \frakg_{\infty} 
= \left\{ if \in \Omega^0(M;i\R) 
  \; \left| \; 
  \mbox{$\int_M if dv_g = 0$} \right\} \right. \\
& \\
& \calT_0|_{\frakg_{\infty}}^* (\cdot) 
= \calT_0^* (\cdot) - \frac{1}{\vol(M)} \int_M \calT_0^*(\cdot) \, dv_g \,. 
\end{array}
\end{equation*} 
Thus the kernel of $\calT_0|_{\frakg_{\infty}}^*$ is the set of those
linearised configurations, which are mapped to $i\R$ under $\calT_0^*$:
\begin{equation*}
\ker \calT_0|_{\frakg_{\infty}}^* 
= \left\{ \left( \left. \begin{array}{c} \nu \\ \phi \end{array} \right) 
  \; \right| \; 
  \calT_0^* \left( \, 
  \left( \begin{array}{c} \nu \\ \phi \end{array} \right) \, \right) \in i\R
  \right\}
= (\calT_0^*)^{-1}(i\R) \; .
\end{equation*}
As above we derive from the complex
$\calK^{\infty}_{(A,\psi)}$ the following $L^2$-orthogonal splitting:  
\begin{align}
\Omega^1(M;i\R) \times \Gamma(\Sigma^+)
&= (\ker \calT_0|_{\frakg_{\infty}}^* \cap \ker \calT_1)  
\oplus \im  \calT_0|_{\frakg_{\infty}} \oplus \im  \calT_1^* \nonumber \\
&= ((\calT_0^*)^{-1}(i\R) \cap \ker \calT_1) 
   \oplus \calT_0(\frakg_{\infty}) \oplus \im \calT_1^* \label{ellsplitP}\; .  
\end{align} 
To define the quotient $L^2$-metric on $\frakP$, we identify the tangent space
\begin{equation*}
T_{[A,\psi]}\frakP 
= T_{[A,\psi]}\left(\preM/\calG_{\infty}\right)
= \ker \calT_1 / \calT_0(\frakg_{\infty})
\end{equation*}
via the splitting splitting \eqref{ellsplitP} with the orthogonal complement
of $\calT_0(\frakg_{\infty})$ in $\ker \calT_1$:
\begin{equation*}
T_{[A,\psi]}\frakP 
\cong \big( (\calT_0(\frakg_{\infty}))^{\perp} \subset \ker \calT_1 \big)
= \ker \calT_0|_{\frakg_{\infty}}^* \cap \ker \calT_1
= (\calT_0^*)^{-1}(i\R) \cap \ker \calT_1 \; .
\end{equation*}
This identification together with the orthogonal splitting \eqref{ellsplitP}
defines a natural Riemannian metric on 
$\preM / \calG_{\infty} \cong \frakP$, which will be called the {\it quotient
  $L^2$-metric on $\frakP$}.

\begin{Lem}
The bundle projection $\frakP \rightarrow \M_{\mu^+}$ is a Riemannian
submersion with respect to the quotient $L^2$-metrics.
\end{Lem}
\begin{proof}
We need to identify the tangent space in $[A,\psi]_{\infty}$ of the fibre
$\frakP_{[A,\psi]}$ over $[A,\psi]$ inside the tangent space 
$T_{[A,\psi]_{\infty}}\frakP$. Then we need to show that
$T_{[A,\psi]}\frakP$ splits orthogonally into the tangent space of  
the fibre $\frakP_{[A,\psi]}$ over $[A,\psi]$ and the tangent space of
$\M_{\mu^+}$, and that the linearisation of the bundle projection 
$\pi: T_{[A,\psi]_{\infty}}\frakP \rightarrow T_{[A,\psi]}\M_{\mu^+}$
is the orthogonal projection of that splitting.

Since $U(1)$ acts on the bundle 
$\preM / \calG_{\infty} \rightarrow \M_{\mu^+}$ via the standard gauge action
of  $U(1) \subset \calG$ on $\preM$, the tangent space of the $U(1)$-orbit
through  $[A,\psi]_{\infty} \in \preM / \calG_{\infty}$  is the image of
$\calT_0(i\R)$ under the quotient map 
$\ker \calT_1 \rightarrow \ker \calT_1 / \calT_0(\frakg_{\infty})$. 
Thus in our model 
$\ker \calT_1 / \calT_0(\frakg_{\infty}) 
\cong \ker \calT_0|_{\frakg_{\infty}}^* \cap \ker \calT_1$,
the tangent space of the fibre $\frakP_{[A,\psi]}$ over
$[A,\psi]$ is the image of $\calT_0(i\R)$ under the orthogonal
projection 
\mbox{$\pi^{\perp}: \ker \calT_1 \rightarrow 
\ker \calT_0|_{\frakg_{\infty}}^* \cap \ker \calT_1$.}
Since $\ker \calT_0^* \subset \calT_0|_{\frakg_{\infty}}^*$, the projection
$\pi^{\perp}$ is the identity on
\begin{equation*} 
T_{[A,\psi]}\M_{\mu^+}
\cong \ker \calT_0^* \cap \ker \calT_1
\; \subset \; 
\ker \calT_0|_{\frakg_{\infty}}^* \cap \ker \calT_1 
\cong T_{[A,\psi]_{\infty}}\frakP \; .
\end{equation*}
Since $\calT_0(i\R)$ is orthogonal to $\ker \calT_0^* \cap \calT_1$, its image
under $\pi^{\perp}$ stays orthogonal to $\ker \calT_0^* \cap \ker \calT_1$. 
Thus the tangent space $T_{[A,\psi]_{\infty}}\frakP_{[A,\psi]}$ of the fibre
over $[A,\psi]$ can be identified with the orthogonal complement of 
$T_{[A,\psi]}\M_{\mu^+}\cong \ker \calT_0^* \cap \ker \calT_1$ in the tangent
space of the total space
$T_{[A,\psi]_{\infty}}\frakP 
\cong \ker \calT_0|_{\frakg_{\infty}}^* \cap \ker \calT_1$. 
This orthogonal complement can be made explicit using the $0$-th order Green
operator $G_0$ of the elliptic complex $\calK_{(A,\psi)}$. Namely, the image
of $\calT_0$ splits $L^2$-orthogonally as 
$\im \calT_0 = \calT_0 \circ G_0 (i\R) \oplus \calT_0 (\frakg_{\infty})$, and
consequently the orthogonal complement of $\ker \calT_0^* \cap \calT_1$ in
$(\calT_0)^{-1}(i\R) \cap \ker \calT_1$ is $\calT_0 \circ G_0 (i\R)$. 
We thus obtain the following $L^2$-orthogonal splitting:
\begin{eqnarray*}
\Omega^1(M;i\R) \times \Gamma(\Sigma^+)
&=& ((\calT_0^*)^{-1}(i\R) \cap \ker \calT_1)
   \oplus \calT_0(\frakg_{\infty}) 
   \oplus \im \calT_1^* \\
&=& \underbrace{
     \underbrace{(\ker \calT_0^* \cap \ker \calT_1)}_{\cong \, T_{[A,\psi]}\M} 
   \oplus  
     \underbrace{\calT_0 \circ G_0 (i\R)}_{\cong \,
        T_{[A,\psi]_{\infty}}\frakP_{[A,\psi]}} 
   }_{\cong \, T_{[A,\psi]_{\infty}}\frakP}
   \oplus \calT_0(\frakg_{\infty})
   \oplus \im \calT_1^* \; .
\end{eqnarray*}
Thus the tangent space of $\frakP$ splits $L^2$-orthogonally as:
\begin{equation*}
T_{[A,\psi]_{\infty}}\frakP
= (\ker \calT_0^* \cap \ker \calT_1) 
\oplus \calT_0 \circ G_0 (i\R) 
\cong T_{[A,\psi]}\M_{\mu^+}
\oplus T_{[A,\psi]_{\infty}}\frakP_{[A,\psi]} \; .
\end{equation*}
It follows that the restriction to $\ker \calT_0^* \cap \ker \calT_1$ of the
linearisation of the bundle projection  $\frakP \rightarrow \M_{\mu^+}$ is an
isometry onto $T_{[A,\psi]}\M_{\mu^+}$. Hence the bundle projection is a
Riemannian submersion as claimed. 
\end{proof}

From the identification of the tangent space
$T_{[A,\psi]_{\infty}}\frakP_{[A,\psi]}$ of the fibre 
over $[A,\psi]$ with {$\calT_0 \circ G_0(i\R)$}, we deduce that the
fundamental vector field $\widetilde{X} \in \frakX(\frakP)$
induced by an element {$X \in \Lie U(1) = i\R$} is given by: 
\begin{equation*}
\widetilde{X}_{[A,\psi]_{\infty}} 
= \big(\calT_0 \circ G_0\big)_{(A,\psi)} (X) \; ,
\end{equation*}
where the subscript indicates the dependence of the operators $\calT_0$ and
$G_0$ on the monopole $(A,\psi) \in [A,\psi]_{\infty}$.

\subsection{The curvature of the quotient $L^2$-metric on $\calB^*$}
\label{chap2.4}
Following the outline above, we compute an explicit formula
for the sectional curvature of the quotient $L^2$-metric on the space
$\calB^*$ of gauge equivalence classes of irreducible configurations in
terms of the Green operators of the elliptic complex $\calK_{(A,\psi)}$ using
the O'Neill formula for the Riemannian submersion $\calC^* \rightarrow \calB^*$. 
It expresses the sectional curvature of the target space $\calB^*$ as the sectional curvature of the source $\calC^*$ plus a positive correction term in the commutator of horizontal extensions $\overline{X}, \overline{Y} \in \frakX (\calC^*)$ of vector fields 
$X,Y \in \frakX (\calB^*)$.
The tangent space in $[A,\psi]$ of the quotient
$\calB^* = \calC^* / \calG$ can naturally be identified with the {horizontal}
space $\calH_{(A,\psi)} = \ker \calT_0^*$. 
Tangent vectors $X_0,Y_0 \in \calH_{(A,\psi)}$, represented as linearised configurations by
\begin{equation*}
X_0 = \left( \begin{array}{c} \nu^X \\ \phi^X \end{array} \right)
\qquad \mbox{resp.} \qquad
Y_0 = \left( \begin{array}{c} \nu^Y \\ \phi^Y \end{array} \right)
\; \in \Omega^1(M;i\R) \times \Gamma(\Sigma^+) \; ,
\end{equation*}
are extended to horizontal vector fields $\overline{X},\overline{Y}$ on $\calC$ simply by projecting the constant extension to the horizontal subbundle:
\begin{equation}
\overline{X}_{(A,\psi)} := \hor_{(A,\psi)} (X_0)
\qquad \mbox{resp.} \qquad
\overline{Y}_{(A,\psi)} := \hor_{(A,\psi)} (Y_0) \; .
\end{equation}
As the proof of the O'Neill formula (see e.g.~\cite{cheeger-ebin75}) relies
only on the algebraic properties of the curvature (such as the Koszul formula)
and on the submersion properties, the formula holds true even in the infinite
dimensional case of the Riemannian submersion $\calC^* \rightarrow \calB^*$.
For $X_0,Y_0 \in T_{[A,\psi]}\calB^* = \ker \calT_0^*$, the O'Neill formula 
reads: 
\begin{equation} \label{oneill}
\big( R^{\calB^*} (X,Y)Y,X \big)_{[A,\psi]} \,
= \, \big( R^{\calC^*} (\overline{X},\overline{Y})\overline{Y},\overline{X}
     \big)_{(A,\psi)} \,
+ \, \frac{3}{4} \,\big\|
     \verti_{(A,\psi)} [\overline{X},\overline{Y}]_{(A,\psi)} \big\|^2 \; .
\end{equation}
The terms of this formula can be computed using the formulae \eqref{vert-hor} for the
orthogonal projectors $\verti_{(A,\psi)}$ and $\hor_{(A,\psi)}$. Since the Levi-Civita connection $\nabla$ of the $L^2$-metric is torsionfree, we have:
\begin{equation*}
[\overline{X},\overline{Y}]_{(A,\psi)}
= \nabla_{X_0} \overline{Y} - \nabla_{Y_0} \overline{X} \; .
\end{equation*}
We may thus express the commutator term in \eqref{oneill}
by covariant derivatives:
\begin{align*}
\big( \nabla_{X_0} \overline{Y} \big)_{(A,\psi)} 
&= \frac{d}{dt} \Big|_0 \overline{Y}\big(\, (A,\psi) + t \cdot X_0 \,\big) \\
&= \frac{d}{dt} \Big|_0 \mbox{hor}_{(A,\psi) + t \cdot X_0} (Y_0)          \\
&= \frac{d}{dt} \Big|_0 \Big( Y 
   - \big\{ \calT_0(t) \circ G_0(t) \circ \calT_0^*(t) \big\} Y_0 \Big) \; , \\
\intertext{where the variable $t$ indicates, that the linearisations and Green
  operators are taken in the point {$(A,\psi) + t \cdot X_0$}. At the initial
  point $t=0$, we just write $\calT_0$ etc. instead of $\calT_0 (t=0)$. Recall
  that $Y_0$ was supposed to be tangent to $\calB^*$, 
  i.e.~$Y_0 \in \ker \calT_0^*$. Hence using the product rule we need only
  differentiate the operator next to $Y_0$, and we thus get:} 
&= - \frac{d}{dt} \Big|_0 \big\{ 
   \calT_0 \circ G_0 \circ \calT_0^*(t) \big\} (Y_0) \\
&= - \frac{d}{dt} \Big|_0 \big\{ \calT_0 \circ G_0 \big\} 
   \big( 2 d^* \nu^Y
   + i \Ima \langle \psi+t\cdot\phi^X,\phi^Y \rangle \big) \\
&= - \calT_0 \circ G_0 \; i \Ima \langle \phi^X,\phi^Y \rangle \; .
\end{align*}
Consequently, the commutator reads: 
\begin{equation*}
[ \overline{X},\overline{Y} ]_{(A,\psi)} 
= - 2 \calT_0 \circ G_0 \; i \Ima \langle \phi^X,\phi^Y \rangle.
\end{equation*}
This term is already vertical, since the vertical bundle
is $\calV = \im \calT_0$. Recall that the $L^2$-metric on $\calC^*$ is flat, so
that we find for the sectional curvature of the space of 
equivalence classes of irreducible connections:
\begin{eqnarray}
\big( R^{\calB^*} (X,Y)X,Y \big)_{[A,\psi]} \,
&=& \frac{3}{4} \, \big\|  
 - 2 \calT_0 \circ G_0 \, i \Ima \langle \phi^X,\phi^Y \rangle \big\|^2_{L^2}
 \notag \\
&=& 3 \, \big( \, 
   \calT_0 \circ G_0 \, i \Ima \langle \phi^X,\phi^Y \rangle \, , \, 
   \calT_0 \circ G_0 \, i \Ima \langle \phi^X,\phi^Y \rangle \, \big)_{L^2} 
 \notag \\
&=& 3 \, \big( \, 
   \calT_0^* \circ \calT_0 \circ G_0 \,
   i \Ima \langle \phi^X,\phi^Y \rangle \, , \,
   G_0 \, i \Ima \langle \phi^X,\phi^Y \rangle \, \big)_{L^2} 
 \notag \\
&=& 3 \big( \, i \Ima \langle \phi ^X,\phi^Y \rangle  \, , \, 
   G_0 \, i \Ima \langle \phi^X,\phi^Y \rangle \, \big)_{L^2} \; .
\end{eqnarray}

\subsection{The curvature of the $L^2$-metric on the
  \mbox{premoduli} space} \label{chap2.5}
The Gauss equation expresses the sectional curvature of a submanifold
with the induced metric in terms of the sectional curvature of the
ambient space and the second fundamental form of the embedding. The
proof of the Gauss equation relies only on the algebraic properties of the
Riemannian curvature tensor and on the definitions of the induced Levi-Civita
connection on the submanifold and of the second fundamental form. Those can
easily be defined in our case using the orthogonal projections onto the
tangent resp.~normal space of the submanifolds discussed above. 
Thus the Gauss equation holds true even for the
$L^2$-metric on the embedding $\preM \hookrightarrow \calC^*$ 
resp.~$\M_{\mu^+}\hookrightarrow \calB^*$.
For the premoduli space $\preM$ the Gauss equation reads:
\begin{eqnarray}\label{gausspremod}
\big( R^{\preM}(X,Y)Y,X \big)_{(A,\psi)} 
&=& \big( R^{\calC^*} 
 (\overline{X},\overline{Y})\overline{Y},\overline{X} \big)_{(A,\psi)}
 \notag \\ 
&& - \big( \II (X,X),\II (Y,Y) \big)_{(A,\psi)}
   + \big( \II (X,Y),\II (X,Y) \big)_{(A,\psi)} \; ,
\end{eqnarray}
where the second fundamental form is defined as 
$\II (X,Y)_{(A,\psi)}:= \big( \nor_{(A,\psi)} (\nabla_X \overline{Y}) \big)$. 
In order to compute the terms of \eqref{gausspremod}, we start with tangent
vectors $X_0, Y_0 \in T_{(A,\psi)}\preM = \ker \calT_1$, represented as
linearised configurations by
\begin{equation*}
X_0 = \left( \begin{array}{c} \nu^X \\ \phi^X \end{array} \right) 
\qquad \mbox{resp.} \qquad
Y_0 = \left( \begin{array}{c} \nu^Y \\ \phi^Y \end{array} \right) 
\; \in \Omega^1(M;i\R) \times \Gamma(\Sigma^+) 
\end{equation*}
and locally extend them to vector fields 
$\overline{X}, \overline{Y}$ on $\calC$ via:
\begin{equation*}
\overline{X}_{(A,\psi)} := \tan_{(A,\psi)} (X_0) 
\qquad \mbox{resp.} \qquad
\overline{Y}_{(A,\psi)} := \tan_{(A,\psi)} (Y_0) \; .
\end{equation*}
Note that $\overline{X}, \overline{Y}$ are indeed extensions to $\calC^*$ of
vector fields on $\preM$: namely, when $(A,\psi)$ is a monopole, then
$\overline{X}_{(A,\psi)} \in T_{(A,\psi)} \preM$. For the covariant derivative
$\nabla_{X_0} \overline{Y}$ we find: 
\begin{align*}
\big( \nabla_{X_0} \overline{Y} \big)_{(A,\psi)} 
&= \frac{d}{dt} \Big|_0 \overline{Y}\big(\, (A,\psi) + t \cdot X_0 \,\big) \\ 
&= \frac{d}{dt} \Big|_0 
   \tan_{(A,\psi) + t \cdot X_0} (Y_0)   \\ 
&= \frac{d}{dt} \Big|_0 \Big( 
   Y_0 - \big\{ \calT_1^*(t) \circ G_2(t) \circ \calT_1(t) \big\} Y_0 \Big)
   \;.\\ 
\intertext{\noindent Recall that $Y_0$ was supposed to be tangent to $\preM$, 
  i.e.~$Y_0 \in \ker \calT_1$. Hence using the product rule we need only
  differentiate the operator next to $Y_0$, and we thus get:}   
&= - \frac{d}{dt} \Big|_0 \big\{ 
   \calT_1^* \circ G_2 \circ \calT_1(t) \big\}Y_0 \\
&= - \frac{d}{dt} \Big|_0  \calT_1^* \circ G_2 \, 
   \left( \begin{array}{c} d^+\nu^Y - q(\psi + t\phi^X,\phi^Y) \\
   \frac{1}{2} \nu^Y \cdot (\psi + t\phi^X) + \dirac_{A + t\nu^X} \phi^Y 
   \end{array} \right) \\
&= - \, \calT_1^* \circ G_2 \, 
   \left( \begin{array}{c} - q (\phi^X,\phi^Y) \\
   \frac{1}{2} \nu^Y \cdot \phi^X + \frac{1}{2} \nu^X \cdot \phi^Y 
\end{array} \right) \; .
\end{align*}
This term is already normal, since the normal space in $(A,\psi)$ is
$N_{(A,\psi)}\preM = \im \calT_1^*$. We thus get for the second fundamental
form terms in the Gauss equation: 
\begin{eqnarray*}
\big( \, \II (X,X),\II (Y,Y) \big)_{(A,\psi)} 
&=&
\left( 
  \calT_1^* \circ G_2 
    \left( \begin{array}{c} - q (\phi^X,\phi^X) \\
    \nu^X \cdot \phi^X \end{array} \right) \; , \;
  \calT_1^* \circ G_2 
    \left( \begin{array}{c} - q (\phi^Y,\phi^Y) \\
    \nu^Y \cdot \phi^Y \end{array} \right)
\right)_{L^2}  \\
& & \\
&=&
\left( 
  \calT_1 \circ \calT_1^* \circ G_2 
    \left( \begin{array}{c} - q (\phi^X,\phi^X) \\
    \nu^X \cdot \phi^X  \end{array} \right) \; , \;
   G_2 
     \left( \begin{array}{c} - q (\phi^Y,\phi^Y \\
     \nu^Y \cdot \phi^Y \end{array} \right)
\right)_{L^2}  \\
& & \\ 
&=&
\left(  
  {}
    \left( \begin{array}{c} - q (\phi^X,\phi^X) \\
    \nu^X \cdot \phi^X \end{array} \right)  \; , \;
   G_2 
     \left( \begin{array}{c} - q (\phi^Y,\phi^Y) \\
     \nu^Y \cdot \phi^Y \end{array} \right) 
\right)_{L^2} \\
\mbox{and analogously:}\qquad \qquad & & \\
& & \\
\big( \, \II (X,Y),\II (X,Y) \big)_{(A,\psi)} &=& 
\left(
    \left( \begin{array}{c} - q (\phi^X,\phi^Y) \\
    \frac{1}{2} \nu^X \cdot \phi^Y + \frac{1}{2} \nu^Y \cdot \phi^X 
    \end{array} \right) \; , \;
  G_2 
    \left( \begin{array}{c} - q (\phi^X,\phi^Y) \\
    \frac{1}{2} \nu^Y \cdot \phi^X + \frac{1}{2} \nu^X \cdot \phi^Y 
    \end{array} \right)
\right)_{L^2} 
\end{eqnarray*}
Since the $L^2$-metric on the configuration space $\calC$ is flat,
we find for the sectional curvature of the premoduli space the
formula:
\begin{eqnarray*}
\big( R^{\preM}(X,Y)Y,X \big)_{(A,\psi)} 
&=& {}  - \big( \II (X,X),\II (Y,Y) \big)_{(A,\psi)}
 + \big( \II (X,Y),\II (X,Y) \big)_{(A,\psi)} \notag \\
&=&{} - \left( 
    \left( \begin{array}{c} - q (\phi^X,\phi^X) \\
    \nu^X \cdot \phi^X \end{array} \right) \; , \; 
  G_2 
    \left( \begin{array}{c} - q (\phi^Y,\phi^Y) \\
    \nu^Y \cdot \phi^Y \end{array} \right)
\right)_{L^2} \notag \\
&&{} + \left( 
    \left( \begin{array}{c} - q (\phi^X,\phi^Y) \\
    \frac{1}{2} \nu^Y \cdot \phi^X + \frac{1}{2} \nu^X \cdot \phi^Y 
    \end{array} \right) \; , \;
  G_2 
    \left( \begin{array}{c} - q (\phi^X,\phi^Y \\
    \frac{1}{2} \nu^Y \cdot \phi^X + \frac{1}{2} \nu^X \cdot \phi^Y 
    \end{array} \right) 
\right)_{L^2}  
\end{eqnarray*}

\subsection{The curvature of the quotient $L^2$-metric on the moduli space} 
\label{chap2.6}
Since the Levi-Civita connections on the quotients $\calB^*$ resp. $\M_{\mu^+}$ are expressed in terms of orthogonal projections from the Levi-Civita connection on $\calC^*$, we can do the computations of the terms in the sectional curvature of $\M_{\mu^+}$ on the configuration space $\calC^*$.
To this end, we consider the family of vector spaces
\begin{equation*}
\calE = \ker (\calT_0^* \oplus \calT_1) 
= \ker \calT_0^* \cap \ker \calT_1 
\rightarrow \calC^* \; .
\end{equation*}
The restriction of $\calE$ to the premoduli space $\preM$ gives a vector bundle of rank $d=-\chi(\calK_{(A,\psi)})$, naturally isomorphic to the pullback of the tangent bundle
of the moduli space: 
$$
\calE|_{\preM} \cong \pi^* T\M_{\mu^+}\; .
$$ 
Note that the dimension of $\ker \calT_0^* \oplus \calT_1$ is not necessarily
constant, hence in general, $\calE$ does not define a vector bundle neither on
the whole configuration space $\C$ nor on the its irreducible part $\calC^*$.
However, the operator operator $\calT_0^* \oplus \calT_1$ is elliptic for every configuration $(A,\psi) \in \calC$, and its index $d=\chi(\calK_{(A,\psi)})$ is independent of $(A,\psi)$. The index equals the dimension of $\ker \big( \calT_0^* \oplus \calT_1 \big)$ minus the dimensions of the obstruction spaces. 
$\calE$ thus defines a vector bundle on the set of those
configurations, for which the obstruction spaces vanish. 

Two tangent vectors $X_0,Y_0 \in T_{[A_0,\psi_0]}\M_{\mu^+}$, represented as linearised configurations by
\begin{equation*}
X_0 = \left( \begin{array}{c} \nu^X \\ \phi^X \end{array} \right) 
\qquad \mbox{resp.} \qquad
Y_0 = \left( \begin{array}{c} \nu^Y \\ \phi^Y \end{array} \right)
\; \in \Omega^1(M;i\R) \times \Gamma(\Sigma^+) \; ,
\end{equation*}
are extended to sections $\overline{X}, \overline{Y}$ of
$\calE$ as: 
\begin{equation} \label{defextension}
\overline{X}_{(A,\psi)} 
:= \tan_{(A,\psi)} \circ \hor_{(A,\psi)} (X_0) 
\qquad \mbox{resp.} \qquad
\overline{Y}_{(A,\psi)} 
:= \tan_{(A,\psi)} \circ \hor_{(A,\psi)} (Y_0) \; .
\end{equation}
We could also have chosen the orthogonal projectors $\tan_{(A,\psi)}$ and
$\hor_{(A,\psi)}$ in reversed order. 
Thus we should keep in mind whether our formulae depend on the choice of the extension. 

Now we proceed as above for the premoduli space to compute the terms
of the Gauss equation. For the covariant derivatives 
$\nabla_{X_0} \overline{Y}$ we find: 
\begin{align} 
\big( \nabla_{X_0} \overline{Y} \big)_{[A,\psi]} 
&= \frac{d}{dt} \Big|_0
   \overline{Y}\big( \, (A,\psi) + t \cdot X_0 \, \big) \nonumber \\
&= \frac{d}{dt} \Big|_0 \tan_{(A,\psi) + t \cdot X_0} \, 
   \circ \, \hor_{(A,\psi) + t \cdot X_0} \, (Y_0) \nonumber \\
&= \frac{d}{dt} \Big|_0 \Big( Y_0 \,
   - \, \big\{ \calT_1^*(t) \circ G_2(t) \circ \calT_1(t) \big\} \, Y_0 \,
   - \, \big\{ \calT_0(t) \circ G_0(t) \circ \calT_0^*(t) \big\} \, Y_0 
   \nonumber \\
&\quad \qquad {} + \big\{ 
   \calT_1^*(t) \circ G_2(t) \circ \calT_1(t) \circ \calT_0(t) \circ G_0(t)
   \circ \calT_0^*(t) 
   \big\} \, Y_0 \, \Big) \nonumber 
\intertext{Recall that $Y_0$ was supposed to be tangent to $\M_{\mu^+}$,
  i.e.~$Y_0 \in 
  \ker \calT_0^* \cap \ker \calT_1$. Hence using the procuct rule, we need
  only differentiate the operators next to $Y_0$, and we thus get:}
&= \frac{d}{dt} \Big|_0 \Big( Y_0 \,
   - \, \big\{ \calT_1^* \circ G_2 \circ \calT_1(t) \big\} \, Y_0 \,
   - \, \big\{ \calT_0 \circ G_0 \circ \calT_0^*(t) \big\} \, Y_0 \nonumber \\
&\quad \qquad {} + \big\{ 
   \calT_1^* \circ G_2 \circ \calT_1 \circ \calT_0 \circ G_0 \circ \calT_0^*(t)
   \big\} \, Y_0 \, \Big) \nonumber 
\intertext{Since $\calK_{(A,\psi)}$ is a complex, we have
  $\calT_1 \circ \calT_0 \equiv 0$, thus the last term vanishes
  identically. If we would have chosen the operators $\tan_{(A,\psi)}$ and
  $\hor_{(A,\psi)}$ in reversed order in \eqref{defextension}, then we would
  have got the term $\frac{d}{dt} \Big|_0 \big\{ \calT_0 \circ G_0 \circ
  \calT_0^* \circ \calT_1^* \circ G_2 \circ \calT_1(t) \big\} \, Y_0$
  instead. But this vanishes by the same argument, since  
  $\calT_0^* \circ \calT_1^* \equiv 0$. We thus get:}
&= \frac{d}{dt} \Big|_0 \Big( {}
  - \big\{ \calT_1^* \circ G_2 \circ \calT_1(t) \big\} \, Y_0 \,  
  - \, \big\{ \calT_0 \circ G_0 \circ \calT_0^*(t) \big\} \, Y_0  \Big) 
\label{nablaXY} \; .
\end{align}
For the second fundamental form terms we need to take the normal projection 
$\nor_{(A,\psi)}$ thereof. Since 
$\big( \im \calT_0 \subset \ker \calT_1 \big) \perp \im \calT_1^*$, the last
term of \eqref{nablaXY} vanishes under $\nor_{(A,\psi)}$
whereas the first term of \eqref{nablaXY} -- being already normal -- stays
unaffected. We thus get for the second fundamental form of the embedding 
$\M_{\mu^+}\hookrightarrow \calB^*$:
\begin{equation}
\II(X,Y)_{[A,\psi]} 
= - \, T_1^* \circ G_2 \, 
  \left( \begin{array}{c} - q (\phi^X,\phi^Y) \\
  \frac{1}{2} \nu^Y \cdot \phi^X + \frac{1}{2} \nu^X \cdot \phi^Y 
  \end{array} \right) \; .
\end{equation}
Note that, although the formulae for the second fundamental forms of the embedding 
$\preM \hookrightarrow \calC^*$ resp.~$\M_{\mu^+} \hookrightarrow \calB^*$ 
look exactly the same, the linearised configurations $X_0, Y_0$ in these formulae are not the same but lie in the different subspaces $\ker \calT_1$ 
resp.~$\ker \calT_0^* \cap \ker \calT_1$ of $T_{(A,\psi)} \calC^*$. 

To proceed we need only collect the terms of the Gauss equation as
for $\preM$ above and combine them with the formula \eqref{oneill} for the
sectional curvature of $\calB^*$. We finally get the following
formula for the sectional curvature of the Seiberg-Witten moduli space
with respect to the quotient $L^2$-metric:
\begin{eqnarray*}
\big( R^{\M}(X,Y)Y,X \big)_{[A,\psi]} 
&=& \big( R^{\calB^*} 
(\overline{X},\overline{Y})\overline{Y},\overline{X} \big)_{[A,\psi]} \\
&& - \big( \II(X,X),\II(Y,Y) \big)_{[A,\psi]}
 + \big( \II(X,Y),\II(X,Y) \big)_{[A,\psi]} \nonumber \\
&=& \, 3 \, \big( \, i \Ima \langle \phi ^X,\phi^Y \rangle  \, , \, 
 G_0 i \Ima \langle \phi^X,\phi^Y \rangle \, \big)_{L^2} \nonumber \\
&& 
- \left( 
    \left( \begin{array}{c} - q (\phi^X,\phi^X) \\
    \nu^X \cdot \phi^X \end{array} \right) \; , \; 
 G_2
    \left( \begin{array}{c} - q (\phi^Y,\phi^Y) \\
    \nu^Y \cdot \phi^Y \end{array} \right) 
\right)_{L^2} \nonumber \\
&&
+ \left(
    \left( \begin{array}{c} - q (\phi^X,\phi^Y) \\
    \frac{1}{2} \nu^Y \cdot \phi^X + \frac{1}{2} \nu ^X \cdot \phi^Y 
    \end{array} \right) \; , \;
  G_2 
    \left( \begin{array}{c} - q (\phi^X,\phi^Y) \\
    \frac{1}{2} \nu^Y \cdot \phi^X + \frac{1}{2} \nu^X \cdot \phi^Y 
    \end{array} \right)
\right)_{L^2}  \; .
\end{eqnarray*}
Note that all these formulae for the sectional curvature
implicitly depend on the perturbation $\mu^+ \in \Omega^2_+(M;i\R)$ used in
the construction of the moduli space. This dependence is via the monopoles
$(A,\psi)$, where our computations are based. These monopoles clearly change,
when the perturbation $\mu^+$ changes.

As the nonlocal Green operators cannot be computed explicitly, we are not
able to draw any direct consequences out of formulae of this type. The best
one can hope for, is that some regularisation techniques allow to compute
e.g.~regularised traces of these operators or that one can compute the terms
more explicitly in special situations. The same problem arises in Yang-Mills
theory, where {\sc Maeda, Rosenberg} and {\sc Tondeur} used 
regularised traces to study the geometry of the gauge orbits in
\cite{maeda-rosenberg-tondeur93,maeda-rosenberg-tondeur95,maeda-rosenberg-tondeur97}, 
whereas {\sc Groisser} and {\sc Parker} used the identification of the
Yang-Mills moduli space of $G=SU(2)$ on $\sph^4$ with instanton number $1$
with the hyperbolic 5-space to compute the curvature of the $L^2$-metric in
the standard instanton $A_0$ explicitly, 
see \cite{groisser-parker87,groisser-parker89}. 
They found, that the $L^2$-metric is {\it not} the standard hyperbolic metric,
but that the curvature in $A_0$ is $\frac{5}{16 \pi^2} > 0$.    

\subsection{The quotient $L^2$-metric on the parametrised moduli space}
\label{chap2.7}
In this section, we construct a natural $L^2$-metric on the parametrised moduli
space $\widehat{\M}$ in the same way as we did for the
moduli space $\M_{\mu^+}$, via appropriate $L^2$-orthogonal splittings.
We show that the restriction of this quotient $L^2$-metric on $\widehat{\M}$
to a fibre $\M_{\mu^+(t_0)}$ of the parametrisation 
$\widehat{\M} = \bigsqcup_{t \in [0,1]} \M_{\mu^+(t)}$ coincides with the
quotient $L^2$-metric of $\M_{\mu^+(t_0)}$, at least if $\M_{t_0}$ is a smooth
manifold. 

The parametrised moduli space $\widehat{\M}$ was defined as the
disjoint union of the moduli spaces $\M_{\mu^+(t)}$ along a curve 
$[0,1] \rightarrow \Omega^2_+(M;i\R), t \mapsto \mu^+(t)$. 
For a generic choice of the curve $t \mapsto \mu^+$, the space $\widehat{\M}$ is a smooth manifold.  
We may further assume that for every $t \in [0,1]$, the derivative $(\mu^+)'_t$ of the curve 
$t \mapsto \mu^+(t)$ is nontrivial. i.~e.~$(\mu^+)'_t \neq 0$. 
We then consider $\widehat{\M}$ as the quotient by $\calG$ of the zero
locus of the parametrised Seiberg-Witten map
\begin{equation*}
\begin{array}{cccc}
\widehat{\calSW}_{\mu^+}:
& \conn \times \Gamma(\Sigma^+) \times [0,1]
& \rightarrow
& \Omega^2_+(M;i\R) \times  \Gamma(\Sigma^-) \\
& & &  \\
& \left( \begin{array}{c} A \\ \psi \\ t \end{array} \right)
& \mapsto
& \left( \begin{array}{c} F_A^+ - \frac{1}{2}q(\psi,\psi) - \mu^+(t) \\ 
  \dirac_A \psi \end{array} \right)  \;.
\end{array}
\end{equation*}
Here, the gauge group $\calG$ acts trivially on the $[0,1]$-factor of 
$\calC \times [0,1]$. 
When we linearise $\widehat{\calSW}_{\mu^+}$ and the orbit map, we end up with the following
complex:
\begin{equation} \tag*{$\widehat{\calK}_{(A,\psi,t_0)}$}  
0 \longrightarrow
  \Omega^0(M;i\R)
  \stackrel{\widehat{\calT}_0}{\longrightarrow}
  \Omega^1(M;i\R) \times \Gamma(\Sigma^+) \times \R
  \stackrel{\widehat{\calT}_1}{\longrightarrow}
  \Omega^2_+(M;i\R) \times \Gamma(\Sigma^-)
  \longrightarrow 
0 
\end{equation} 
Here $\widehat{\calT_1}$
denotes the linearisation in $(A,\psi,t_0)$ of the parametrised Seiberg-Witten
map $\widehat{\calSW}_{\mu^+}$
\begin{equation*}
\begin{array}{cccc}
\widehat{\calT}_1:
& \Omega^1(M;i\R) \times \Gamma(\Sigma^+) \times \R
& \rightarrow
& \Omega^2_+(M;i\R) \times \Gamma(\Sigma^-) \\
& & & \\
& \left( \begin{array}{c} \nu \\ \phi \\ s \end{array} \right)
&\mapsto
& \left( \begin{array}{c}
  d^+\nu - q(\psi,\phi) - s \cdot (\mu^+)'_{t_0} \\
  \frac{1}{2} \nu \cdot\psi + \dirac_A \phi \end{array} \right) 
\end{array}
\end{equation*}
and $\widehat{\calT}_0$
denotes the linearisation in ${\bf 1} \in \calG$ of the orbit map through
$(A,\psi,t_0)$
\begin{equation*}
\begin{array}{cccc}
\widehat{\calT}_0:
& \Omega^0(M;i\R)
& \rightarrow
& \Omega^1(M;i\R) \times \Gamma(\Sigma^+) \times \R \\
& & & \\
& if
& \mapsto
& \left( \begin{array}{c} 2idf \\  -if\cdot\psi \\ 0  \end{array} \right) \;.
\end{array}
\end{equation*}
Note that $\widehat{\calK}_{(A,\psi,t_0)}$ is a complex, since
$\widehat{\calT}_1 \circ \widehat{\calT}_0 \equiv \calT_1 \circ \calT_0 \equiv
0$ holds trivially, but it is {\it not} elliptic. 
However, the splittings used in section \ref{chap2.2} to construct the quotient
$L^2$-metric can still be obtained directly from these operators and their
adjoints. 
The adjoint of $\widehat{\calT}_0$ is the operator: 
\begin{equation} \label{hatT0*}
\begin{array}{cccc}
\widehat{\calT}_0^*:
& \Omega^1(M;i\R) \times \Gamma(\Sigma^+) \times \R
& \rightarrow
& \Omega^0(M;i\R) \\
& & & \\
& \left( \begin{array}{c} \nu \\ \phi \\ s \end{array} \right)
& \mapsto
& 2d^*\nu \, + \, i \Ima \langle \psi,\phi \rangle 
\end{array}
\end{equation} 
and the adjoint of
$\widehat{\calT}_1$
is the operator:
\begin{equation} \label{hatT1*}
\begin{array}{cccc}
\widehat{\calT}_1^*:
& \Omega^2_+(M;i\R) \times \Gamma(\Sigma^-)
& \rightarrow
& \Omega^1(M;i\R) \times \Gamma(\Sigma^+) \times \R \\
& & & \\
& \left( \begin{array}{c} \mu \\ \xi \\ s \end{array} \right)
& \mapsto
& \left( \begin{array}{c}
  d^*\mu + \frac{i}{2} \Ima \langle (\cdot) \cdot \psi,\xi \rangle \\
  \dirac_A \xi - 2 \mu \cdot \psi \\ - s \cdot ((\mu^+)'_{t_0},\mu)_{L^2} 
  \end{array} \right) \;.
\end{array}
\end{equation}
Although the complex $\widehat{\calK}_{(A,\psi)}$ is not elliptic, the
operators $\widehat{\calT}_0$ and $\widehat{\calT}_1^*$ are obviously closed,
and we thus have the following $L^2$-orthogonal splitting:
\begin{eqnarray} 
\Omega^1(M;i\R) \times \Gamma(\Sigma^+) \times \R
&=& \ker \widehat{\calT}_0^* \oplus \im \widehat{\calT}_0 \notag \\
&=& \ker \widehat{\calT}_1 \oplus \im \widehat{\calT}_1^* \notag \\
&=& (\ker \widehat{\calT}_0^* \cap \widehat{\calT}_1) 
  \oplus \im \widehat{\calT}_0 \oplus \im \widehat{\calT}_1^* \label{parasplit}
\end{eqnarray}
Similar to the case of the moduli space as explained in section \ref{chap2.2},
the intersection of the kernels of $\widehat{\calT}_0^*$ and
$\widehat{\calT}_1$ can be regarded as the Zariski tangent space of
the parametrised moduli space $\widehat{\M}$. Thus in an irreducible point
$\widehat{[A,\psi]}$, the parametrised moduli space $\widehat{\M}$ carries a
natural Riemannian metric, induced from the splitting \eqref{parasplit} and the
identification 
$T_{\widehat{[A,\psi]}}\widehat{\M} 
\cong \ker \widehat{\calT}_0^* \cap \ker \widehat{\calT}_1$. 
As before, we call the metric obtained in this way the {\it quotient}
$L^2$-{\it metric}
on the parametrised moduli space
$\widehat{\M}$. 

To show that the $L^2$-metric induced from the embedding 
$\M_{\mu^+_t} \hookrightarrow \widehat{\M}$ coincides with $L^2$-metric on the moduli space $\M_{\mu^+(t)}$ as constructed in section \ref{chap2.2}, we compare the images and kernels of $\widehat{\calT}_j^{(*)}, j = 0,1$ with those of $\calT_j^{(*)}, j = 0,1$, and we find:
\begin{eqnarray*}
\ker \calT_0^* \times \R
&=& \ker \widehat{\calT}_0^*  \\
\im \calT_0 \times \{0\} 
&=& \im \widehat{\calT}_0 
 \\
\ker \calT_1 \times \{0\}
&=& \ker \widehat{\calT}_1 
  \, \cap \, 
  \Omega^1(M;i\R) \times \Gamma(\Sigma^+) \times \{0\} \\
\im \calT_1^* \times \R 
&=& \im \widehat{\calT}_1^* \,.
\end{eqnarray*}
The tangent space to a regular slice 
$\M_{\mu^+_{t_0}} \subset \widehat{M} = \bigsqcup_{t \in [0,1]} \M_{\mu^+(t)}$
can thus be identified with the intersection of the tangent space of $\widehat{\M}$ with the tangent space of the $t_0$-slice in $\calC^* \times [0,1]$: 
\begin{eqnarray*}
T_{[A,\psi]}\M_{\mu^+(t_0)} 
&=& \ker \calT_0 \cap \ker \calT_1^* \\
&\cong& (\ker \widehat{\calT}_0 \cap \widehat{\calT}_1^*) 
\, \cap \, \Omega^1(M;i\R) \times \Gamma(\Sigma^+) \times \{0\} \\
&=& T_{\widehat{[A,\psi]}}\widehat{\M} 
\, \cap \, \Omega^1(M;i\R) \times \Gamma(\Sigma^+) \times \{0\} \; .
\end{eqnarray*}
The restriction of the quotient $L^2$-metric on $\widehat{\M}$ to a regular slice 
$\widehat{\M}_{t_0} = \M_{\mu^+_{t_0}}$ thus yields the natural
quotient $L^2$-metric constructed in section \ref{chap2.2}.

Finally extracting the material from the subsections, we can summarise our
results in the theorem:

\begin{Thm}
Let $M$ be a compact smooth 4-manifold with a fixed \SpinC{} and
$\mu^+$ (resp.~$\mu^+(t), t \in [0,1]$) generic perturbations such that the
Seiberg-Witten moduli space $\M_{\mu^+}$ (resp.~the parametrised moduli space
$\widehat{\M} = \bigsqcup_{t \in [0,1]} \M_{\mu^+(t)}$) are smooth manifolds
of the expected dimension. 
Then there exists a natural quotient $L^2$-metric on $\M_{\mu^+}$ and a
compatible quotient $L^2$-metric on $\widehat{\M}$ such that the metric
induced from the inclusion of a smooth slice $\M_{\mu^+(t_0)} \hookrightarrow
\widehat{\M}$ is the same as the metric constructed on $\M_{\mu^+(t_0)}$ as
the moduli space with perturbation $\mu^+(t_0)$. 
In case $M$ is simply connected, the Seiberg-Witten bundle $\frakP \to
\M_{\mu^+}$ -- i.e.~the isomorphism class of principal $U(1)$ bundles on
$\M_{\mu^+}$ defining the invariants -- admits a
natural geometric representative carrying a quotient $L^2$-metric such that
the projection $\frakP \to \M_{\mu^+}$ is a Riemannian submersion.
The sectional curvature of those metrics is explicitly given in terms of
the Green operators of the deformation complex of the Seiberg-Witten
equations.
\end{Thm}

The above construction automatically yields a natural $L^2$-metric on the
regular part $\M^*_{\mu^+}$ of the (perturbed or nonperturbed) moduli
spaces. Hence if one does not want to bother with the problem of how to choose
appropriate perturbations in order that the smoothness obstructions vanish and
the moduli spaces be regular, one could restrict to the regular part to obtain
a Riemannian metric on $\M^*_{\mu^+}$. However, there is no reason to hope for
this metric to be complete. 


\section{Moduli spaces on K\"ahler surfaces} \label{chap3}
In this section, we recall the well known identification of Seiberg-Witten monopoles on K\"ahler surfaces with vortices.
We further recall the identification of the moduli space $\M_{\mu^+}$ as a torus fibration over the complex projective space $\Pro( H^0_{A_0}(M;L))$ and as a K\"ahler quotient of a K\"ahler submanifold of the configuration space, which follows from the work of {\sc Bradlow} and {\sc Garc{\'\i}a-Prada} on the vortex equations on compact K\"ahler manifolds in \cite{bradlow90,bradlow91,garcia-prada_1994a,garcia-prada98}.
That the Seiberg-Witten moduli spaces appear as symplectic quotients had also been remarked by {\sc Okonek} and {\sc Teleman} in \cite{okonek-teleman95a}.
A far more general statement for all kinds of moduli spaces, which range under the universal
Kobayashi-Hitchin correspondance, has been established by {\sc L\"ubke} and
{\sc Teleman} in \cite{luebke_teleman_06}.

\subsection{Seiberg-Witten equations on K\"ahler surfaces}\label{chap3.1}
On a K\"ahler surface, the (perturbed) Seiberg-Witten equations take a very
simple form in terms of holomorphic data. As pointed out by {\sc Witten} in
\cite{witten94}, these are special kinds of so called vortex equations, which
had been first studied by {\sc Bradlow} in \cite{bradlow90,bradlow91} and by
{\sc Garc{\'\i}a-Prada} in \cite{garcia-prada93,garcia-prada_1994b,garcia-prada_1994a}. 
For a detailed discussion of the relation between the Seiberg-Witten and vortex equations,
see also \cite{bradlow_garcia-prada_1995,garcia-prada98}. 
Equivalently, the monopoles can be identified in terms of algebraic geometry
as effective divisors. 
From an existence theorem for solutions of the vortex equations, {\sc Witten}
first deduced the nontriviality of the invariants on a K\"ahler surface. 
However, this identification of monopoles with vortices resp.~effective
divisors holds some subtleties: it depends on the choice of a rather special
type of perturbations. 
In that case, no Sard-Smale argument is available to
ensure the regularity of the moduli spaces for a generic perturbation. 
Thus it may happen, that the moduli spaces are not generically smooth, or that they
are smooth manifolds but not of the expected dimension $d= - \chi(\calK_{(A,\psi)})$. 
However, the obstructions for the moduli space to be smooth of the expected dimension near a monopole can be identified rather explicitly in terms of cohomology groups associated with the
effective disivor corresponding to that monopole, see
\cite{friedman_morgan_1999} and \cite{okonek-teleman96}.   

Since this hardly affects our consideration of $L^2$-metrics, we do not recall those
obstructions here.
Instead, we restrict the study of the geometry to the
regular part $\M_{\mu^+}^*$ of the moduli space. 
To assure that the obstruction spaces vanish and the moduli spaces be regular, one could also restrict the consideration to more special K\"ahler surfaces, such
as ruled surfaces with additional properties, see \cite{friedman_morgan_1999}. 

Throughout this section let $(M,g)$ be a compact, connected K\"ahler surface
with K\"ahler form $\omega$. The complex structure determines a canonical
\SpinC{} $P_0$, whose determinant line bundle is the dual of the
canonical line bundle $K_M = \Lambda^{2,0}T^*M$,
i.e.~$\det P_0= K_M^* = \Lambda^{0,2}T^*M$.
Any other \SpinC{} $P$ has the form $P = P_0 \otimes L$ for a $U(1)$-bundle $L$, and the determinant line bundle of $P$ is then given by $\det P= K_M^* \otimes L^2$.
We will not distinguish in notation between a $U(1)$-bundle
$L$ and its associated complex line bundle.  
The positive resp.~negative spinor bundles are:
\begin{equation*}
\Gamma(\Sigma^+) = \Omega^0(M;L) \oplus \Omega^{0,2}(M;L) 
\qquad \mbox{and} \qquad 
\Gamma(\Sigma^-) = \Omega^{0,1}(M;L) \; .
\end{equation*}
Let $A_L$ be a connection on the line bundle $L$ and $A_{can}$ the Chern connection, i.e.~the unique hermitean holomorphic connection on $\Lambda_M^*$. 
The Dirac operator of the \SpinC{} $P = P_0 \otimes L$ with respect to the product connection 
$A = A_{can} \otimes A_L^2$ on $\det (P)$ is given by
$\dirac_A 
= \sqrt{2} ( \overline{\partial}_{A_L} + \overline{\partial}^*_{A_L} )$.

Taking constant multiples $\mu^+ = i\pi\lambda \cdot \omega$, $\lambda \in \R$ of the K\"ahler form as perturbations (these are clearly transversal to the wall $\Gamma^+_g$), the Seiberg-Witten equations read: 
\begin{eqnarray}
\big( F_A^+ \big)^{1,1}
&=&  \frac{i}{4} \big( |\beta|^2 - |\zeta|^2 \big) 
      \cdot \omega \, + \, i \pi \lambda \omega  \label{SWK1}  \\
\big( F_A^+ \big)^{0,2}
&=&  \frac{ \overline{\beta} \zeta }{2} \label{SWK2}  \\
\sqrt{2} \big( \overline{\partial}_{A_L}\beta 
   + \overline{\partial}^*_{A_L}\zeta \big)
&=&  0 \label{SWK3}
\end{eqnarray}
As shown by {\sc Witten} in \cite{witten94}, for a monopole 
$(A,\beta \oplus \zeta), \beta \in \Omega^0(M;L), \zeta \in \Omega^{0,2}(M;L)$, one of the components $\alpha,\zeta$ necessarily vanishes.
Which one vanishes, is detected by the degree
$\deg_{\omega}(L) := \int_M c_1(L) \wedge \omega$
of the line bundle $L$.
The corresponding result for the equations perturbed (or ``twisted'', as they insist) by a closed real $(1,1)$-form $\mu^+$ has been established by {\sc Okonoek} and {\sc Teleman} in \cite{okonek-teleman95a}.
Namely, for perturbations $\mu^+ = i\pi\lambda \cdot \omega$, $\lambda \in \R$, we have $\zeta \equiv 0$, if $\lambda \cdot \vol(M) \leq \deg_\omega(\det P)$ and $\alpha \equiv 0$, if 
$\lambda \cdot \vol(M) \geq \deg_\omega(\det P)$.
In either case, the determinant line bundle $\det P= K_{M}^* \otimes L^2$ 
carries the structure of a holomorphic line bundle, and with
respect to the induced holomorphic structure on $L$, the
components  $\beta$ resp.~$\overline{\zeta}$ are holomorphic sections of $L$
resp.~$K_{M} \otimes L^*$.
By replacing the line bundle $L$ with $K_M \otimes L^*$ if neceessary,
one can always arrange $\deg_{\omega} (\det P)$ to
have a fixed sign. 
Hence by choosing $\lambda \leq \frac{\deg_{\omega} (\det P)}{\vol(M)}$, we may assume the monopoles to be of the form $(A,\beta) \in \calA(\det P) \times \Omega^0(M;L)$.

{\sc Witten} further observed in \cite{witten94}, that for K\"ahler surfaces
with $b_2^+ >1$, when taking holomorphic 2-forms as perturbations, the only
generically nonempty moduli spaces are those of the canonical and the
anticanonical \SpinC{}. Consequently, for all but those two \SpinCs{}, the
Seiberg-Witten invariant vanishes. 
The situation is completely different in the case $b_2^+ = 1$: as shown by {\sc Okonek} and {\sc Teleman} in \cite{okonek-teleman96}, the Seiberg-Witten invariants of a K\"ahler manifold $M$ with $b_2^+(M) = 1$ and $b_1(M) = 0$ are nontrivial in precisely one chamber, as soon as the \SpinC{}
has nonnegative index.

\subsection{Monopoles and vortices}
We briefly recall the identification of Seiberg-Witten monopoles (for
perturbations $\mu^+ = i \pi \lambda \omega, \lambda \in \R$) with vortices
resp.~effective divisors, as first established by {\sc Witten}, and later made precise by {\sc Okonek} and {\sc Teleman} in \cite{okonek-teleman95a,okonek-teleman96} and by {\sc Friedman} and {\sc Morgan} in \cite{friedman_morgan_1999}. 
This identification yields an isomorphism of real analytic spaces between 
Seiberg-Witten moduli spaces and Douady spaces of effective divisors of fixed 
topological type, see e.g.~\cite{okonek-teleman96}. 
In the regular case, this isomorphism is a diffeomorphism to the complex projective space 
$|D_0| = \Pro( H^0(M;[D_0]))$ for $b_1(M) = 0$ resp. to a fibration through complex projective spaces over the torus $H^1(M;i\R)/H^1(M;2 \pi i\Z)$ of holomorphic structures.

A connection $A \in \calA(L)$ is called {\it holomorphic}, if the $(0,2)$-part $F_A^{0,2}$ of its curvature vanishes.
Any two holomorphic connections differ by a $1$-form of type $(1,0)$, thus the space of holomorphic connections is an affine space modelled over the complex vector space $\Omega^{1,0}(M)$, see e.g.~\cite{huybrechts_2005}.
A holomorphic connection on a line bundle $L$ determines a holomorphic structure via the Cauchy-Riemann operator
\begin{equation*}
\overline{\partial}_A := \frac{1}{2} (d_A - i \cdot d_A \circ \frakJ) \,.
\end{equation*}
The {\it complexified gauge group}
$\calG^{\C} := \Omega^0(M;\C^*)$ acts on $\calA(L)$ by
\begin{equation*}
u: A \mapsto A 
+ u^{-1}\overline{\partial}u - \bar{u}^{-1}\partial\bar{u} \;,
\end{equation*}
which extends the action of the gauge group $\calG$. 
The induced action on the Cauchy-Riemann operator by 
\begin{equation*}
u: \overline{\partial}_A \mapsto \overline{\partial}_{u \cdot A}
= u^{-1} \circ \overline{\partial} \circ u \; .
\end{equation*}
yields isomorphisms of the complex structures.
The space of isomorphism classes of holomorphic structures on $L$ may thus be regarded
as the quotient of the space of holomorphic connections
$\calA^{hol}(L)$ by the action $\calG^{\C}$. 
This quotient can be identified with the torus $H^1(M;i\R)/H^1(M;2\pi i\Z)$, see
e.g.~\cite{salamon95}.

The complexified gauge group acts on the sections of $L$ by $\alpha \mapsto u^{-1} \cdot \alpha$. 
Obviously, the section $\alpha \in \Omega^0(M;L)$ is holomorphic with respect to
$\overline{\partial}_A$ if and only if the section $u^{-1} \cdot \alpha$ is
holomorphic with respect to $\overline{\partial}_{u \cdot A}$.  
Any holomorphic connection $A$ on $L$ is related by a complex gauge
transformation $u$ to the hermitean holomorphic (or Chern-) connection $A_0$
of the corresponding holomorphic structure. 
Writing a complex gauge transformation $u \in \calG^{\C}$ as
$u = e^{-f+ih}$ with real functions $f,h$ we find: 
\begin{eqnarray*}
A = u \cdot A_0 
&=& A_0 + u ^{-1} \overline{\partial} u - \bar{u}^{-1} \partial \bar{u} \\
&=& A_0 + \overline{\partial} (-f +ih) - \partial (-f -ih) \\
&=& A_0 + id^c\!f + idh \; .
\end{eqnarray*}

A Seiberg-Witten monopole consists of a holomorphic connection $B$ on $L$ and a $B$-holomorphic section $\beta$, which satisfy \eqref{SWK1}.
Via a complex gauge transformation $u$, we may write $(B,\beta)$ as
\begin{equation*}
  \left( \begin{array}{c}
  B \\ \beta \end{array} \right) 
= u \cdot \left( \begin{array}{c} 
  A_0 \\ \alpha \end{array} \right)
= \left( \begin{array}{c} 
  A_0 + 2id^c\!f + 2idh \\ e^{f-ih} \alpha \end{array} \right)
\end{equation*}
with $\alpha \in H^0_{A_0}(M;L)$ -- the $A_0$-holomorphic sections -- 
and real functions $f,h$. Since $e^{ih}$ is an ordinary gauge transformation, 
we end up with an equation in $f$:
\begin{equation*}
(F_B^+)^{1,1} 
= F^+_{A_0} + \big( 2idd^c\!f \big)^+  
= \frac{i}{4} |\beta|^2 \cdot \omega + i \pi \lambda \cdot \omega  
= \frac{i}{4} e^{2f} |\alpha|^2 \cdot \omega + i \pi  \lambda \cdot \omega \,.
\end{equation*}
Contracting both sides with the K\"ahler form $\omega$, we see that any Seiberg-Witten monopole 
$(B,\beta) = (A_0 + 2id^c\!f + 2idh, e^{f-ih}\alpha)$ 
can be derived from a configuration
$(A_0,\alpha), \alpha \in H^0_{A_0}(M;L)$ by solving the equation
\begin{equation}
2\Delta f + \frac{1}{2} e^{2f} |\alpha|^2
= - 2 \pi \lambda - i\Lambda_{\omega} (F_{A_0}^+) \; .
\label{laplacef}
\end{equation}
in the unknown $f$. 

In this form, the Seiberg-Witten equations are a special case of so called vortex equations, which had first been studied by {\sc Bradlow} and {\sc Garc{\'\i}a-Prada}.
They gave several proofs for the existence and uniqueness of vortices as well as identifications of the corresponding moduli spaces.
The first proof was given by {\sc Bradlow} in \cite{bradlow90} by using the existence and uniqueness theorem for solutions of equations of the form \eqref{laplacef} due to {\sc Kazdan} and {\sc Warner} in \cite{kazdan-warner74}.
For another proof using the continuity method, see \cite{becker_thesis}.
Different constructions of vortices were given by {\sc Bradlow} and {\sc Garc{\'\i}a-Prada} in \cite{bradlow91,garcia-prada93,garcia-prada_1994a}.

Summarising, for any nonzero $A_0$-holomorphic section $\alpha \in H^0_{A_0}(M;L)$, there is a unique solution $f \in \mathcal{C}^\infty(M)$ to the equation \eqref{laplacef}.
Thus any such section yields a Seiberg-Witten monopole $(B,\beta) = u \cdot (A_0,\alpha)$.
Taking gauge equivalence into account, we obtain for any orbit of $\calG^{\C}$ in $\calA^{hol}(L)$ a bijection from the projective space $\Pro (H^0_{A_0}(M;L))$ to the intersection of the Seiberg-Witten moduli space with that orbit.
In case all regularity obstructions vanish, this bijection is a diffeomorphism and the Seiberg-Witten moduli space thus fibres over the torus 
$H^1(M;i\R)/H^1(M;2\pi i\Z)$ of isomorphism classes of holomorphic structures through the complex projective spaces $\Pro (H^0_{A_0+\nu}(M;L))$, where $A_0$ is the Chern connection of a fixed holomorophic structure and
$[\nu] \in H^1(M;i\R)/H^1(M;2\pi i\Z)$.

\subsection{$b_2^+=1$ and the parametrised moduli space}
As well known from \cite{okonek-teleman96}, the Seiberg-Witten invariants on a
K\"ahler surface $M$ with $b_2^+(M) = 1$ are nonzero in exactly one of the two
chambers, as soon as the \SpinC{} has nonnegative index. If additionally 
$b_1(M) = 0$, then the moduli space for a perturbation $\mu^+$ on the wall
$\Gamma_g^+$ consists of a single point: The Kazdan-Warner type equation for an arbitrary monopole  $(B_1,\beta_1)= (A_0 + 2id^c\!f_1 + 2idh_1 , e^{f_1-ih_1} \alpha_1)$ and a reducible monopole
$(B_2,\beta_2) = (A_0 + 2id^c\!f_1 + 2idh_1 , 0)$ reads:
\begin{equation} \label{KWpoint}
2 \Delta f_1 + \frac{1}{2} e^{2f_1} |\alpha|^2 
= - 2\pi \lambda - i \Lambda_{\omega}(F_{A_0}^+)
= 2 \Delta f_2 \; .
\end{equation}
This implies $2 \Delta(f_1 - f_2) = \frac{1}{2} e^{2f_1} |\alpha_1|^2$. 
By integration, we find 
$\| e^{f_1} \alpha_1 \|^2 \equiv 0$, so $\alpha_1 \equiv 0$. 
Thus $(B_1,\beta_1)$ is a reducible monopole too. 
Equation \eqref{KWpoint} now implies $\Delta (f_1 - f_2) = 0$, 
thus solutions $f_1,f_2$ differ by a constant.
Correspondingly, the monopoles $(B_1,\beta_1)$, $(B_2,\beta_2)$ are gauge equivalent.

We now consider the behaviour of the quotient $L^2$-metric under changes of
the perturbation. When the perturbation $\mu^+$ approaches the wall
$\Gamma^+_g$, the moduli spaces $\M_\mu^+$ collapses from a compact space --
homeomorphic to a complex projective space -- to a point. 
We show that this collaps is indeed a collaps in the quotient $L^2$-metric, i.e. that the
diameter $\diam(\M_{\mu^+})$ of the moduli space in that metric shrinks to
$0$, when the perturbation $\mu^+$ approaches the wall $\Gamma^+_g$.

Choose a path $t \mapsto \mu^+(t)$ of perturbations, such that 
$\mu^+(t) \in \Omega^2_+(M;i\R) - \Gamma^+_g$ for $t \in [t_0,t_1)$ and $\mu^+(t_1) \in \Gamma^+_g$. 
For a generic such path, the moduli spaces $\M_{\mu^+(t)}$ are nonempty and the parametrised moduli space $\widehat{\M} = \bigsqcup_{t \in [t_0,t_1]} \M_{\mu^+(t)}$ is a smooth manifold. 
The topological type of the fibres $\M_{\mu^+(t)}$ collpases from a complex projective space to a point. 
As is well known from the study of wall crossing phenomena (see \cite{okonek-teleman96} and references therein), the parametrised moduli space $\widehat{\M}$ is compact and has the homeomorphism type of a cone on $\M_{\mu^+(t_0)} \cong \CP^m$. 
We denote the tip of the cone, i.e. the unique reducible gauge equivalence class, by $[B',\beta']$. 
The fibre $\M_{t_1} = \{[B',\beta']\}$ may or may not be singular in $\widehat{\M}$.

The quotient $L^2$-metric on the Zariski tangent spaces of the parametrized
moduli space $\widehat{\M}$ is a Riemannian metric on the nonsingular part
$\widehat{\M}^* = \widehat{\M} - \M_{t_1} = \widehat{\M} - \{[B',\beta']\}$.
The Riemannian distance of this metric can be extended in the
singularity $[B',\beta']$, which then has a finite distance from any other
point on $\widehat{\M}^*$. Thus the Riemannian distance makes the parametrised
moduli space into a complete metric space. It is clear, that the diameter of
the fibre $\M_t$ (with respect to this {\it extrinsic} metric) shrinks to $0$,
when the parameter $t$ tends to $t_1$. That the same holds true with respect
to the {\it intrinsic} metric of the fibres, i.e.~the quotient $L^2$-metrics
of $\M_t$, is not a priori clear. Therefor we show:

\begin{Lem}
Let $M$ be a compact K\"ahler surface with $b_1(M) = 0$ and $b_2^+(M) =
1$. Choose a generic path $t \mapsto \mu^+(t)$ of perturbations, such that 
$\mu^+(t) \in \Omega^2_+(M;i\R) - \Gamma^+_g$ for $t \in [t_0,t_1)$ and 
$\mu^+(t_1) \in \Gamma^+_g$ such that the parametrised moduli space 
$\widehat{\M} = \bigsqcup_{t \in [t_0,t_1]} \M_{\mu^+(t)}$ is smooth and of the
expected dimension. Then the diameter $\diam(\M_t)$ of the fibre 
$\M_{\mu^+(t)}$ shrinks to $0$ when the perturbation 
$\mu^+(t)$ approaches the wall $\Gamma^+_g$, i.e. when $t$ tends to $t_1$.
\end{Lem}
\begin{proof}
Suppose this were not the case. Then there would exist an $\epsilon >0$ and a
sequence of points $[B_1,\beta_1]_t, [B_2,\beta_2]_t \in \M_{\mu^+(t)}$ such that
$\dist([B_1,\beta_1]_t,[B_2,\beta_2]_t) = \epsilon \; 
\forall t \in [t_0,t_1)$. 
The points $[B_1,\beta_1]_t, [B_2,\beta_2]_t$ can be joined by geodesics
$\gamma_t$ of length $\epsilon$, and we may take $\gamma_t$ to be parametrised
by arc length. The theorem of Arzela-Ascoli implies that the curves $\gamma_t$
converge uniformly when $t$ tends to $t_1$, and it is clear that the limit
$\gamma_{t_1}$ is the constant curve in $\M_{\mu^+(t_1)} =
\{[B',\beta']\}$. We show that the length of the limit is bounded from below
by $\frac{\epsilon}{2}$: 

The length of the limit curve $\gamma_{t_1}$ in the metric space
$\widehat{\M}$ is defined as the supremum over all
partitions $0= s_0 < \ldots < s_n =\epsilon$ of the interval $[0,\epsilon]$
of the length of the polygon through the points $\gamma_{t_1}(s_i)$:
\begin{equation*} 
\length (\gamma_{t_1}) :=
\sup\limits_{s_0 < \ldots < s_n}
\left( \sum\limits_{i=1}^{n} \dist(\gamma_t(s_{i-1}),\gamma_t(s_i)) \right) \;.
\end{equation*}
For a given partition $s_0 < \ldots < s_n$, we find a parameter 
$t' \in [t_0,t_1]$ sufficiently close to $t_1$ such that
\begin{equation}\label{delta}
\dist(\gamma_t(s_i),\gamma_{t_1}(s_i)) < \delta := \frac{\epsilon}{4n}
\qquad \forall t \in [t',t_1], \forall i = 0 , \ldots , n \; .
\end{equation}
From the triangle inequality and \eqref{delta}, we get:
\begin{equation*}
\dist(\gamma_{t_1}(s_{i-1}),\gamma_{t_1}(s_i))
> \dist(\gamma_t(s_{i-1}),\gamma_t(s_i)) - 2 \delta 
\qquad \forall t \in [t',t_1], \forall i = 0 , \ldots , n \; .
\end{equation*}
We thus obtain the following estimate for the length of the curve
$\gamma_{t_1}$:
\begin{eqnarray*}
\length(\gamma_{t_1}) 
&\geq& \sum\limits_{i=1}^{n} \dist(\gamma_{t_1}(s_{i-1}),\gamma_{t_1}(s_i)) \\
&>& \sum\limits_{i=1}^{n} \dist(\gamma_t(s_{i-i}),\gamma_t(s_i)) - 2\delta\\
&=& \length(\gamma_t) - 2n \delta \\
&=& \epsilon - 2n \cdot \frac{\epsilon}{4n} \\
&=& \frac{\epsilon}{2} \; .
\end{eqnarray*}
This contradicts the fact, that the limit $\gamma_{t_1}$ is the constant curve
in $\M_{\mu^+(t_1)} = \{[B',\beta']\}$ and thus has length $L(\gamma_{t_1}) =
0$.
\end{proof}

\subsection{Moduli spaces as K\"ahler quotients} \label{chap3.4}
In this section we recall the identification of the regular part $\M_{\mu^+}^*$ of the Seiberg-Witten moduli space as a K\"ahler quotient of a certain submanifold of the irreducible configuration space $\calC^*$.
The first of the Seiberg-Witten equations \eqref{SWK1} appears as the zero locus equation of a moment map for the gauge group action on the configuration space $\calC$, whereas the equations \eqref{SWK2}, \eqref{SWK3} define a K\"ahler submanifold $\frakN \subset \calC^*$.
The moment map in question had been computed by {\sc Garc{\'\i}a-Prada} in \cite{garcia-prada_1994a}.
A similar K\"ahler quotient construction appears for moduli spaces of Hermitean-Einstein connections in \cite{kobayashi_1987}.
The relation of vortices to Hermitean-Einstein structures and the corresponding moduli spaces is discussed in \cite{garcia-prada_1994b}.
The conclusion that the K\"ahler quotient construction in the infinite dimensional setting indeed yields a K\"ahler metric on the moduli spaces essentially relies upon the work of {\sc Hitchin} on moduli spaces of vortices resp.~Higgs bundles in \cite{hitchin84,hitchin90a,hitchin90b}. 

The $L^2$-metric on the configuration space $\calC = \conn \times \Omega^0(M;L)$ is a K\"ahler metric with respect to the complex structure
\begin{equation*}
\begin{array}{cccc}
\frakJ^{\calC} 
= \frakJ^{T^*M} \oplus (-i):
&T_{(B,\beta)}\calC = \Omega^1(M;i\R) \times \Omega^0(M;L)
&\rightarrow
&\Omega^1(M;i\R) \times \Omega^0(M;L) \\
& & & \\
&\left( \begin{array}{c} \nu \\ \phi \end{array} \right)
&\mapsto
&\left( \begin{array}{c} \frakJ^{T^*M} \nu \\ 
   (-i) \cdot \phi \end{array} \right) \; .
\end{array}
\end{equation*}
The action of gauge group $\calG$ clearly preserves both the $L^2$-metric and the symplectic form $\Phi^{\calC} = \big( \frakJ^{\calC}\cdot,\cdot \big)_{L^2}$ and has the moment map
\begin{equation*}
\begin{array}{cccc}
\mu^{\calC}:   
&\calC 
&\rightarrow 
&\Omega^0(M,i\R) \subset \frakg^*\\
& & & \\
&\left( \begin{array}{c} B \\ \beta \end{array} \right)
&\mapsto 
&\Lambda_{\omega} (F_B) - \frac{i}{2} |\beta|^2 \,, 
\end{array}
\end{equation*}
were the Lie algebra $\frakg = \Omega^0(M;i\R)$
of the gauge group is identified via the $L^2$-metric as a subset of its dual
$\frakg^*$.
Since $\frakg = \Omega^0(M;i\R)$ is an abelian Lie algebra, we can add any
$\Lambda_\omega(\mu^+)$ with $\mu^+ \in \Omega^2_+(M;i\R)$ to get another moment map 
\begin{equation*}
(B,\beta) 
\mapsto 
\Lambda_{\omega}(F_B) - \frac{i}{2}|\beta|^2 
- \Lambda_{\omega}(\mu^+) \,.
\end{equation*} 
Hence the first Seiberg-Witten equation \eqref{SWK1} appears as the zero locus equation for a moment map on the configuration space. 

The solution space $\frakN$ of the equations \eqref{SWK2} and \eqref{SWK3}
\begin{equation*}
\frakN:= \left\{ \left.
\left( \begin{array}{c} B \\ \beta \end{array} \right) \in \calC \; \right| \; 
\sqrt{2} \, \overline{\partial}_B \beta = 0\,, F_B^{0,2}=0\,,\beta \neq 0 \right\}
\end{equation*}
is a K\"ahler submanifold of the irreducible configuration space $\calC^*$: the restriction of the symplectic form $\Phi^\calC$ to $\frakN$ is nondegenerate and the complex structure $\frakJ^{\calC}$ preserves the tangent bundle 
\begin{equation*}
T\frakN 
= \left\{ \left. 
  \left( \begin{array}{c} \nu \\ \phi \end{array} \right) 
  \in \Omega^1(M;i\R) \times \Gamma(\Sigma^+) \; \right| \; 
  \sqrt{2} \big( \overline{\partial}^*_B \phi + \nu^{0,1} \wedge \beta \big) 
  = 0\,, (d\nu)^{0,2} =0\,, 
  \left( \begin{array}{c} B \\ \beta \end{array} \right) \in \frakN \right\}
\end{equation*}
of $\frakN$. 
The restriction of the moment map $\mu^\calC$ to $\frakN$ gives a moment map $\mu^\frakN$ for the gauge group action on $\frakN$.
The Seiberg-Witten moduli space thus appears as the K\"ahler reduction
\begin{eqnarray*}
\M_{\mu^+} &=& (\mu^{\frakN})^{-1}(0) / \calG  \,.
\end{eqnarray*}
Consequently, the $L^2$-metric on the irreducible configuration space $\calC^*$ descends to a K\"ahler metric on the regular part $\M_{\mu^+}^*$ of the moduli space.
Detailed proofs for the case of moduli spaces of vortices resp.~Higges bundles are due to {\sc Hitchin} \cite{hitchin84,hitchin90a,hitchin90b}.

There is a well known explicit description of the Seiberg-Witten moduli
spaces on the complex projective plane $\CP^2$ with the \SpinC{} 
$P = P_0 \otimes \calO(k)$, $k \in \N$ as 
$$
\M_{\mu^+}(P_0 \otimes \calO(k)) \cong |\calO_{\CP^2}(k)| \cong
\Pro(H^0(M;\calO(k)) \cong \Pro(\C_k[z^0,z^1,z^2])\,,
$$
see \cite{okonek-teleman96}. 
For $k=1$, the quotient $L^2$-metric on the Seiberg-Witten bundle $\frakP \to \M_{\mu^+}$ as constructed in \ref{chap2.2} and \ref{chap2.3} is preserved by the standard $U(3)$-action on $\CP^2$.
Summarising, we thus end up with the following corollary:
\begin{Cor}
Let $M$ be a compact K\"ahler surface. 
Then the quotient $L^2$-metric on the regular part $\M_{\mu^+}^*$ of the Seiberg-Witten moduli space is a K\"ahler metric.
For $M=\CP^2$ with the \SpinC{} $P = P_0 \otimes \calO(1)$, the Seiberg-Witten bundle 
$\frakP \to \M_{\mu^+}$ with the quotient $L^2$-metrics is isometric to the Hopf bundle $S^5 \to \CP^2$ with a Berger metric on $S^5$ and the Fubini-Study metric on $\CP^2$.
\end{Cor}

\bibliographystyle{amsplain}
\bibliography{lit} \addcontentsline{toc}{chapter}{Bibliography}

\end{document}